\documentclass[11pt]{amsart}

\usepackage{amsfonts}
\usepackage{wrapfig}
\usepackage{epstopdf}
\usepackage{enumitem}
\usepackage{algorithmic}

 \usepackage{amsfonts,mathscinet}
 \usepackage{color}
 \usepackage{amsmath, amssymb}
  \usepackage{hyperref}

\usepackage{amsopn}

\usepackage{multirow}
\usepackage{float}
\usepackage[normalem]{ulem}
\usepackage{graphicx}
\usepackage{tikz}
\usepackage{mathrsfs}
\usepackage{extarrows}
\usetikzlibrary{calc}
\usetikzlibrary{matrix,arrows,decorations.pathmorphing}
\usetikzlibrary{arrows,decorations.pathmorphing,decorations.pathreplacing}
\usetikzlibrary{positioning,shapes,shadows,arrows}

\newcommand{\rel}{\protect\tikz[baseline] \protect\draw[line width=.55pt, overlay, line cap=round, dash pattern=on 0 off 1.2pt] (-2.4ex, .4ex) to[bend right=46] (3ex, .4ex);}

\newcommand{\tikztoarg}[1]{\mkern 5.5mu \protect\tikz[baseline]\protect\path[-stealth,line width=.4pt] (0ex,0.65ex) edge node[above=-.4ex, overlay, font=\scriptsize] {$#1$} (3.75ex,.65ex); \mkern 5.5mu}

\numberwithin{equation}{section}

 \newtheorem{theorem}{Theorem}[section]
 \theoremstyle{definition}
 \newtheorem{definition}[theorem]{Definition}

\newcommand{\End}{\operatorname{End}\nolimits}
\newcommand{\thick}{\mathrm{thick}}
\newcommand{\Hom}{\mathrm{Hom}}

\newcommand{\add}{\operatorname{{\rm add }}}

\newcommand{\za}{\alpha}
\newcommand{\zb}{\beta}
\newcommand{\zd}{\delta}
\newcommand{\zD}{\Delta}

\newcommand{\zg}{\gamma}
\newcommand{\zG}{\Gamma}

\newcommand{\Z}{\mathbb{Z}}

\newcommand{\aaa}{{f}}

\newcommand{\proj}{\operatorname{{\rm proj }}}
\newcommand{\tw}{\operatorname{{\rm tw}}}
\newcommand{\Tw}{\operatorname{{\rm Tw}}}

\newcommand{\calt}{\mathcal{T}}

\newcommand{\cala}{\mathcal{A}}

\newcommand{\cald}{\mathcal{D}}

\newcommand{\calw}{\mathcal{W}}

\newcommand{\calo}{\mathcal{O}}

\newcommand{\dba}{\mathcal{D}^b(A)}
\newcommand{\ka}{\operatorname{per}(A)}
\newcommand{\SSS}{\mathbb{S}}
\newcommand{\per}{{\rm per}}
\newcommand{\D}{\mathcal D}

\newcommand{\HH}{\operatorname{HH}}

\definecolor{dark-green}{RGB}{14,150,2}
\definecolor{red}{RGB}{250,0,0}
\newcommand{\gpoint}{\color{black}{\circ}}
\newcommand{\rpoint}{\color{red}{\bullet}}

\begin{document}

\title{On geometric models in  representation theory}
   % \author{Sibylle Schroll }\address{Department of Mathematics, University of Cologne, Cologne, Germany, \email{schroll@math.uni-koeln.de} (\url{https://sites.google.com/site/sibylleschroll/}).}
\author{Sibylle Schroll}
\address{Department f\"ur Mathematik, Universit\"at zu K\"oln, Weyertal 86-90, K\"oln, Germany }
\email{schroll@math.uni-koeln.de}

\date{}

\maketitle

%\fancyfoot[R]{\scriptsize{Copyright \textcopyright\ 20XX by SIAM\\
%Unauthorized reproduction of this article is prohibited}}

\begin{abstract} 
Geometric models have emerged as an important tool in the representation theory of algebras. Surface models associated to gentle algebras have been particularly fruitful in advancing our understanding of their module and derived categories. We give an overview of some of the theoretical advances that geometric surface models for the derived categories of graded gentle algebras and their connections to Fukaya categories of surfaces have made possible.
\end{abstract}

%\tableofcontents

%\maketitle

\section{Introduction}

Triangulated categories and their differential graded  or $\infty$-categorical enhancements provide a unifying framework for algebra, geometry and topology and they have reshaped our understanding of homological algebra. Homological mirror symmetry has been instrumental in building bridges between algebraic and symplectic geometry through derived equivalences between categories of coherent sheaves and Fukaya categories \cite{KontsevichICM}. Tilting theory and exceptional sequences have linked triangulated categories in algebraic and symplectic geometry to the representation theory of algebras \cite{HandBookOfTiltingTheory, Baer, Bondal}.

In the representation theory of associative algebras, triangulated categories appear in many shapes and sizes, for example, as unbounded, bounded or perfect derived categories of modules, as stable module categories or as singularity categories. It is perhaps an unavoidable consequence of their inherent ubiquity and versatility, that the inner workings of derived categories often appear intractable and concrete descriptions are scarce and hard to come by. For derived categories of algebras, this tractability can be quantified by the notions of derived discreteness and derived tameness, meaning that the indecomposable objects in the derived category are,  for every  dimension of the cohomology, either finite in number (derived discreteness) \cite{Vossieck} or they appear in an at most finite number of one-parameter families (derived tameness) \cite{BekkertDrozd, Drozd, GeissKrause}. For derived discrete or derived tame algebras, a classification of the indecomposable objects in their derived category has been achieved for several classes of algebras, notably for gentle and skew-gentle algebras \cite{BekkertMarcosMerklen, BekkertMerklen,  BurbanDrozd, HaidenKatzarkovKontsevich}.

In recent years, a new perspective on derived categories of gentle algebras has greatly advanced our understanding of their derived categories and indeed of triangulated categories in general. Gentle algebras are path algebras of quivers with quadratic monomial relations (see Section \ref{sec: gentle}). The maximal paths in a gentle algebra $A$ naturally appear as vertices of a ribbon graph \cite{Schroll15} whose associated surface $S$ encodes much of the derived representation theory of $A$. Indeed, the perfect derived category $\per (A)$ of $A$ --- or rather its differential graded (dg) enhancement given by the category $\tw (\add A)$ of twisted complexes over  the additive closure  $\add A$ of  $A$  --- is equivalent to the partially wrapped Fukaya category of the surface $S$ \cite{HaidenKatzarkovKontsevich, LekiliPolishchuk} (see also Section \ref{sec: Fukaya}). In general, there is a correspondence 
\begin{equation}
\label{eq:correspondence}
\begin{tikzpicture}[x=1em,y=1em,baseline=-2.25ex]
\node[font=\footnotesize, align=center, anchor=north] (L) at (-12, 0) {\normalsize $A$ \\ graded gentle algebra};
\node[font=\footnotesize, align=center, anchor=north] (R) at (3, 0) {\normalsize $S$ \\ graded marked surface \rlap{ (together with a dissection)}};
\path[<->, line width=.5pt] (-3,-.75) edge (-6,-.75);
\end{tikzpicture}
\end{equation}
  where the grading of $A$ corresponds to a grading of $S$ in the form of  a line field. On the level of objects, this correspondence associates to a curve in $S$ an explicit twisted complex. In particular, the indecomposable direct summands of $A$ can be viewed as curves on the surface, and they function as a coordinate system for general objects of the derived category. By considering infinite twisted complexes, one obtains a model not only for the perfect derived category $\per (A)$, but also for the full bounded derived category $\dba$ \cite{OpperPlamondonSchroll}.

A plethora of additional information about derived categories of gentle algebras can be derived from this geometric perspective. For example, the winding numbers of certain curves on the surface assemble to a complete derived invariant of the category \cite{AmiotPlamondonSchroll,JinSchrollWang, LekiliPolishchuk, Opper}, augmenting a classical combinatorial invariant of gentle algebras \cite{AvellaGeiss}. The action of the Serre functor on objects of $\per(A)$ can be described by changing the endpoints of the corresponding curves on the surface \cite{OpperPlamondonSchroll}, and this description can also be used to compute the categorical entropy of the Serre functor \cite{ChangElaginSchroll}. The Hochschild cohomology, including its full Tamarkin--Tsygan calculus, is likewise encoded by the surface: it corresponds to the boundary components with $0$ or $1$ stops \cite{BarmeierSchrollWangDef, ChaparroSchrollSolotarSuarez}. 

The particulars of the case of derived categories of gentle algebras have also shed light on general phenomena in the study of derived categories. Full exceptional sequences and their associated semiorthogonal decompositions are a key  technology in the structure theory of derived categories. In the derived category of a gentle algebra, exceptional sequences can be characterised by certain dissections of the associated surface \cite{ChangJinSchroll,ChangSchrollexcp}. This geometric perspective has led to the first example of a triangulated category generated by an exceptional sequence such that the natural braid group action on the set of isomorphism classes of full exceptional sequences is not transitive \cite{ChangHaidenSchroll}, answering a longstanding question of Bondal and Polishchuk \cite{BondalPolishchuk}. 

Deformation theory can be used to leave the world of gentle algebras and smooth surfaces. Remarkably, the full algebraic deformation theory of the derived category of a graded gentle algebra, parametrised by its second Hochschild cohomology, can  be understood through its surface model \cite{BarmeierSchrollWangDef}. Hochschild $2$-cocycles correspond to boundary components with $0$ or $1$ stops and winding number $1$ or $2$. Deformations then correspond to partially compactifying some of these boundary components either to smooth points, or to orbifold points, depending on the winding number of the corresponding boundary component. On the one hand, this is the first case confirming the relationship between $A_\infty$-deformations of Fukaya categories and partial compactifications in the context of partially wrapped Fukaya categories. Such a relationship was proposed (in the fully wrapped case) by Seidel in \cite{SeidelICM}. On the other hand, it leads to a theory of Fukaya categories of surfaces with orbifold singularities \cite{AmiotPlamondon,BarmeierSchrollWang,ChoKim}. Just as in the smooth case, these Fukaya categories are categories of global sections of a cosheaf of A$_\infty$-categories \cite{BarmeierSchrollWang}, confirming a singular analogue of a conjecture of Kontsevich \cite{Kontsevich}. Moreover, by characterizing dissections giving rise to formal generators of these categories \cite{AmiotPlamondon,BarmeierSchrollWang,KimPhD}, one obtains a new class of derived tame algebras given by explicit quivers with relations which are derived equivalent to skew-gentle algebras.

\section*{Acknowledgements}
I would like to thank all of my collaborators, irrespective of whether they are cited in this survey or not. Not only have I immensely enjoyed  these collaborations, I also have learned a lot from every single one and this proceedings article is a result of all these experiences. I'm very grateful to have had the opportunity to work with so many outstanding mathematicians.  I also would like to extend a particular thank you to Severin Barmeier who throughout the writing of these proceedings has given unwavering support and feedback.  Furthermore, I would like to thank everyone who has kindly read and given feedback on  earlier versions.

\section{Gentle algebras: Definitions and first properties}\label{sec: gentle}

Growing out of tilting theory and the representation theory of  biserial algebras, dating back to the 1970s and 80s, the literature on gentle algebras and their generalisations is vast and rapidly growing. In this section we briefly recall the definition of gentle algebras and some of their properties.  Unfortunately, it is far beyond the scope of these proceedings to do justice to all the interesting works on this subject.

 \paragraph{Definition of gentle algebras} Gentle algebras were introduced in  \cite{AssemHappel, AssemSkowronski} in the study of iterated tilted algebras of (affine) Dynkin type $\mathbb{{A}}$.

\begin{definition}\label{def:gentle}
{\rm  Let $\Bbbk$ be a field. A  $\Bbbk$-algebra $A$ is \textit{gentle} if it is Morita equivalent to an algebra $\Bbbk Q/I$ for a finite quiver $Q=(Q_0, Q_1)$ such that at every vertex in the set of vertices $Q_0$ of $Q$ there are at most two arrows ending and at most two arrows starting. 
The ideal $I$ is generated by a set of paths of length two such that
\newline
 (1) for any arrow $a$ there is at most one arrow $b$ and at most one arrow $c$ such that $ba \notin I$ and $ac \notin I$,
\newline
 (2) for any arrow $a$ there is at most one (composable) arrow $b'$ and at most one (composable) arrow  $c'$ such that $b'a \in I$  and $a c' \in I$.
\newline
A $\mathbb{Z}$-grading  on $A$ is a  map from the set of arrows $Q_1$ to $\mathbb{Z}$. }
\end{definition}

 Whenever we write algebra, we implicitly understand it to be an algebra over a  base field $\Bbbk$. All statements hold for $\Bbbk$ algebraically closed of characteristic zero, but in many cases, neither of these assumptions is needed.

We remark that a gentle algebra as in Definition~\ref{def:gentle} is not necessarily finite dimensional. 
The following are examples of gentle algebras given by their quivers with relations 
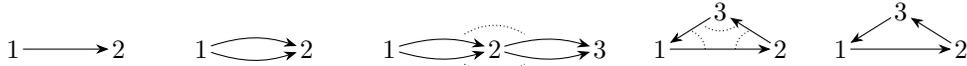
\usetikzlibrary{positioning,arrows.meta}
\[
\begin{tikzpicture}[>=Stealth, baseline=(current bounding box.center)]
\tikzset{every node/.style={font=\small, inner sep=1pt}}

% 1) 1 -> 2
\begin{scope}[shift={(0,0)}]
  \node (a1) at (0,0) {1};
  \node (a2) at (1.4,0) {2};
  \draw[->] (a1) -- (a2);
\end{scope}

% 2) 1 ⇉ 2
\begin{scope}[shift={(2.5,0)}]
  \node (b1) at (0,0) {1};
  \node (b2) at (1.4,0) {2};
  \draw[->,overlay ] (b1) to[bend left=18]  (b2);
  \draw[->] (b1) to[bend right=18] (b2);
\end{scope}

% 3) 1 ⇉ 2 ⇉ 3
\begin{scope}[shift={(5.0,0)}]
  \node (c1) at (0,0) {1};
  \node (c2) at (1.4,0) {2};
  \node (c3) at (2.8,0) {3};
  % double arrows 1 -> 2
  \draw[->,overlay] (c1) to[bend left=16]  (c2);
  \draw[->,overlay] (c1) to[bend right=16] (c2);
  \draw[densely dotted,overlay] (1,0.2) to [bend left] (1.8,0.2);
   \draw[densely dotted,overlay] (1,-0.2) to [bend right] (1.8,-0.2);
  % double arrows 2 -> 3
  \draw[->] (c2) to[bend left=16]  (c3);
  \draw[->] (c2) to[bend right=16] (c3);
\end{scope}

% 4) 3-cycle: 1 -> 2 -> 3 -> 1
\begin{scope}[shift={(8.6,0)}]
  \node (d1) at (0,0)    {1};
  \node (d2) at (1.6,0)  {2};
  \node[overlay] (d3) at (0.8,0.5){3};
  \draw[->] (d1) -- (d2);
  \draw[->,overlay] (d2) -- (d3);
  \draw[->,overlay] (d3) -- (d1);
  \draw[densely dotted,overlay] (0.6,0) to [bend right,overlay] (0.4,0.25);
   \draw[densely dotted,overlay] (1,0) to [bend left,overlay] (1.2,0.25);
     \draw[densely dotted,overlay] (0.6,0.3) to [bend right,overlay] (1,0.3);
\end{scope}

% 5) another 3-cycle: 1 -> 2 -> 3 -> 1
\begin{scope}[shift={(11.0,0)}]
  \node (e1) at (0,0)    {$1$};
  \node (e2) at (1.6,0)  {$2$};
  \node[overlay] (e3) at (0.8,0.5) {$3$};
%  \draw[dotted] 
  \draw[->] (e1) -- (e2);
  \draw[->,overlay] (e2) -- (e3);
  \draw[->,overlay] (e3) -- (e1);
\end{scope}

\end{tikzpicture}
\]
where the dots indicate  zero relations, that is, the corresponding paths of length two are generators of the ideal $I$ in each case. Other examples of gentle algebras are derived discrete algebras (apart from those of Dynkin type $\mathbb{D}$ and $\mathbb{E}$) \cite{Vossieck}, see also \cite{BobinskiGeissSkowronski}, and Jacobian algebras from quivers with potentials from triangulations of marked surfaces with marked points in the boundary \cite{ABCP, Labardini}.

\paragraph{First properties of gentle algebras}

Gentle algebras are quadratic monomial algebras, thus by \cite{GreenZacharia} they are Koszul algebras. It immediately follows from the definition that the Koszul dual, given by the quadratic dual, also is  gentle. Furthermore, gentle algebras are of tame representation type. In \cite{WaldWaschbusch}, see also \cite{ButlerRingel}, building on \cite{GelfandPonomarev, RingelDihedral},  their indecomposable modules are classified combinatorially in terms of strings and bands. These are certain words in the alphabet given by the arrows in the quiver as well as their formal inverses.  An explicit basis of the morphism space between indecomposable modules is given in \cite{CBmorphisms, Krausemorphisms}.   Extension spaces  between indecomposable modules are described in different ways in \cite{ BaurSchroll, BDMTY, CanakciPauksztelloSchrollExtensions, CanakciSchrollExt}. 
In \cite{SchroerZimmermann} it is shown that finite dimensional gentle algebras are closed under derived equivalence, that is, any algebra derived equivalent to a finite dimensional gentle algebra is again gentle. In the catalogue of finite dimensional algebras \cite{SchroerAtlas} this is a rare phenomenon. We note that the closure under derived equivalence is open in the case of graded gentle algebras.

What makes gentle algebras remarkable is that not only their module categories but also their derived categories are tame. The  indecomposable objects in the derived category are again given by a type of word combinatorics, namely they are parametrised by \textit{homotopy strings and bands} \cite{BekkertMerklen, BurbanDrozd, BurbanDrozdNodal1}. Another important class of algebras which is known to be derived tame are skew-gentle algebras \cite{BekkertMarcosMerklen, BurbanDrozd}. Skew-gentle algebras  are, up to Morita equivalence, skew-group algebras of gentle algebras by a $\mathbb{Z}_2$-action. They were originally defined in \cite{GeissdelaPena}.
It is remarkable that, bar a few exceptions, 
all other derived tame algebras that are currently known -- such as the derived tame Nakayama algebras in \cite{BekkertGiraldoVelez}, the derived tame quadratic string algebras in \cite{Fonseca}, and the algebras in \cite[Definition 5.7]{BurbanDrozdnodal} -- are derived equivalent to skew-gentle algebras  \cite{KimPhD}, see also \cite{AmiotPlamondon, BarmeierSchrollWang, ChoKim}. 

In general, derived categories are difficult to parse and explicit calculations are often impossible. However, for the derived category of gentle algebras this is different. For the case of derived discrete gentle algebras a complete treatment of their derived categories is given in \cite{BroomheadPauksztelloPloogII,BroomheadPauksztelloPloogI}. It is also in this context that a first geometric model in terms of the representation theory is formulated in \cite{Broomhead} giving rise to the classification of thick subcategories of the bounded derived category of derived discrete algebras. For the more general case of finite dimensional gentle algebras of finite global dimension, a classification of the thick subcategories generated by string objects of the bounded derived category   has recently been given in \cite{Page} in terms of the geometric model in \cite{OpperPlamondonSchroll}. The classification of thick subcategories generated by band objects is still open. 

The explicit description of bases of the morphism spaces in the bounded derived category of a gentle algebra in \cite{ArnesenLakingPauksztello} was an important step towards our understanding of these derived categories. In \cite{Bobinski}, using the Happel functor from the bounded derived category to the stable module category of the repetitive algebra,  the Auslander--Reiten triangles are described. In \cite{ArnesenLakingPauksztello}, a description of the irreducible morphisms in terms of the basis morphisms are given.  
Another important concept to understand in a triangulated category is that of a mapping cone. For gentle algebras the mapping cones of the basis morphisms in \cite{ArnesenLakingPauksztello} are  calculated in \cite{CanakciPauksztelloSchrollMappingCone}.

\section{Geometric surface models for gentle algebras} \label{sec: geometric model intro}

In this section we give a brief overview of different types of geometric models that can be associated to gentle algebras. This includes  models coming from cluster theory, models for the module category and a model for the derived category.  

\paragraph{Geometric surface models and cluster theory} One important avenue in which geometric surface models  came to the representation theory of algebras is  through their connection with cluster algebras and their categorifications, see, for example, \cite{FominShapiroThurston}. 
In this setting of surface cluster algebras, through the work of \cite{ABCP, CalderoChapotonSchiffler, Labardini} it quickly became apparent that gentle algebras play an important role. They appear as Jacobian algebras of  quivers with potential coming from  ideal triangulations of  surfaces with marked points in the boundary. In \cite{ABCP} it is shown, that the surface encodes the representation theory of the corresponding gentle algebra. In particular, a bijection between homotopy classes of curves connecting marked points (excluding the arcs in the triangulation) and string modules over the gentle algebra is established as well as a correspondence of  closed curves with 1-parameter families of band modules. A similar result is shown in \cite{BruestleZhang} for cluster categories arising from surfaces with marked points in the boundary using the representation theory of gentle algebras. For cluster categories from surfaces with punctures an analogous result is proved in \cite{QiuZhou17,Schiffler} using the representation theory of skew-gentle algebras. In particular,  in these models, extensions between indecomposable objects correspond to crossings of the associated curves, see also \cite{CanakciSchrollExt}. Cluster categories from triangulations of surfaces arise as Verdier quotients of perfect derived categories of 3-Calabi-Yau Ginzburg dg algebras associated to quivers with potential obtained from triangulations.
For these derived categories there exist several geometric models such as the decorated marked surfaces studied,  for example, in \cite{IkedaQiu, Qiu16, QiuZhou19}, and via perverse schobers in \cite{Christ}.

\paragraph{Geometric surface models for module categories} Extending the surface model in \cite{ABCP}, a model for the module category of an arbitrary gentle algebra is given in \cite{BaurCoelho}, with a generalisation to skew-gentle algebras given in \cite{HeZhouZhu}. Using  these models, the $\tau$-tilting theory of gentle and skew-gentle algebras is studied in \cite{ HeZhouZhu,PaluPilaudPlamondonMem}. In \cite{Chang}, the algebraic hearts of the bounded derived categories of gentle algebras are studied in terms of their geometric models given in \cite{OpperPlamondonSchroll}. This gives, in  particular, a different model for the module category which is compatible with that of the bounded derived category described in Section \ref{sec: OPS Model}.

\paragraph{Towards a geometric surface model for the derived category: Ribbon graphs and surface dissections} Coming back to a geometric model for the derived category of   gentle algebras, 
it is shown in \cite{AssemSkowronski,  PogorzalySkowronski, RingelRepetitive}, see also \cite{Schroer}, that a finite dimensional algebra $A$ is gentle if and only if its trivial extension is special biserial. The class of special biserial algebras is a large class of finite dimensional algebras of tame representation type \cite{WaldWaschbusch}. Furthermore, trivial extensions of finite dimensional algebras are symmetric algebras, where an   algebra $A$ is symmetric if there is an isomorphism of $A$-$A$-bimodules from $A$ to its $\Bbbk$-linear dual $DA$. By \cite{AntipovGeneralov, Roggenkamp, Schroll15} every symmetric special biserial algebra is a  so-called \textit{Brauer graph algebra}. This 
%chain of reasoning 
shows that $A$ is gentle if and only if its trivial extension is a Brauer graph algebra. A Brauer graph algebra is 
%an algebra 
defined on a (decorated) ribbon graph $\zD$. In \cite{MarshSchroll}, see also \cite{SchrollBrauer}, a compact oriented surface with boundary and a surface dissection into polygons is constructed from the ribbon graph $\zD$ corresponding to the Brauer graph of  a Brauer graph algebra. It is shown that flips of diagonals of an underlying triangulation or $m$-angulation correspond to derived equivalences of the corresponding Brauer graph algebras. A full derived invariant for Brauer graph algebras in terms of the surface datum is given in \cite{OpperZvonareva} using a generalisation of the covering theory in \cite{GreenSchrollSnashall} and an approach via  $A_\infty$-algebras.  For a more general construction in this direction, see \cite{ChristHaidenQiu} and for a generalisation of gentle and Brauer graph algebras, see \cite{GreenSchrollBrauerconfig, GreenSchrollalmostgentle}.

It follows from the above, that one can associate to every gentle algebra a compact oriented surface via the ribbon graph $\zD$ which is the Brauer graph of the Brauer graph algebra given by the trivial extension. It is then a natural question  whether the ribbon graph of this Brauer graph algebra can directly be read from the gentle algebra. That this is the case is shown in \cite{Schroll15}, where the ribbon graph $\zD$  is defined directly from the gentle algebra without passing to the trivial extension. It is precisely this ribbon graph $\zD$ which in \cite{OpperPlamondonSchroll}  gives rise to  the representation theoretic formulation of a geometric model for the bounded derived category of a gentle algebra in terms of a surface dissection into polygons. This surface dissection is obtained by gluing the vertices of the ribbon graph $\zD$ to the boundary of the corresponding surface. In this way $\zD$ can be viewed as  a surface dissection into polygons. 
In the context of partially wrapped Fukaya categories, $\zD$ corresponds  to a full formal arc system defined in \cite{HaidenKatzarkovKontsevich}. We note that a dual version of dissecting a surface into polygons is given in \cite{CoelhoParsons} under the name of tiling algebra. The tiling of the surface can be seen as a dual ribbon graph (dual to $\zD$) of a gentle algebra and corresponds to the (dual) surface dissection $\zD^*$ in Section \ref{sec: OPS Model}. This dual ribbon graph of a gentle algebra is also constructed in  \cite{LekiliPolishchuk} in the context of partially wrapped Fukaya categories of surfaces with stops, see Section \ref{sec: Fukaya}. It is shown in \cite{OpperPlamondonSchroll}  that these two dual constructions related to $\zD$ and $\zD^*$ are linked by  Koszul duality, see also \cite{LiQiuZhou} for a functorial description. In summary, the different constructions above lead to a bijection between surface dissections into polygons up to homeomorphism of surfaces and gentle algebras up to isomorphism \cite{BaurCoelho, OpperPlamondonSchroll, PaluPilaudPlamondon}.

\section{Derived categories of graded gentle algebras and their geometric model}\label{sec: OPS Model} 

In this section we briefly introduce the geometric model of the bounded derived category of a graded gentle algebra in \cite{OpperPlamondonSchroll}. We note that \cite{OpperPlamondonSchroll} was first developed for ungraded finite dimensional gentle algebras (of not necessarily finite global dimension). In a later version the results were extended to  the  bounded derived category $\dba$ of a finite dimensional not necessarily homologically smooth graded gentle algebra $A$, considered as a proper dg algebra with zero differential, using, in particular, unbounded twisted complexes recently studied in \cite{AL}. 

Below we describe more explicitly the correspondence in \eqref{eq:correspondence} as summarised in the following table 

{ \small
\begin{center}
\begin{tabular}{ c|c } 
 \hline
 \rule{0pt}{2.2ex} $\dba$, $A$ proper graded gentle  & Compact oriented surface with $\gpoint$-points in the boundary  \\[.1em]
  \hline\hline \rule{0pt}{2.2ex}
 \smash{\raisebox{-6pt}{$X_\zg$ indecomposable object}}  & $\zg$  graded arc or \\
 &  graded closed curve with an indecomposable local system \\[.1em]
 \hline \rule{0pt}{2.2ex}
 $\Bbbk$-basis of $\Hom_{\dba}(X_\zg, X_\zd)$ & oriented intersections from $\zg$ to $\zd$  \\[.1em]
 \hline \rule{0pt}{2.2ex}
 \smash{$\raisebox{-6pt}{$X_\zg \stackrel{f}\to X_\zd \to X_{\za_1} \oplus X_{\za_2}  \to$}$} & resolution of the crossing of $\zg$ and $\zd$ corresponding to $f$ resulting in \\
 & two curves $\za_1$ and $\za_2$ (one of which might be zero) \\[.1em]
 \hline \rule{0pt}{2.2ex}
\smash{\raisebox{-6pt}{Serre functor $\SSS$ on $\ka$}}
& moving both endpoints of a graded arc connecting $\gpoint$-points to \\ & the next $\gpoint$-point on 
 the respective boundary components  \\[.2em]
\hline
\end{tabular}
\end{center}
}

\subsection{Bounded and perfect derived categories}\label{sec: twisted complexes}

Before we describe the geometric model and the different correspondences in more detail, we set up the categorical framework that we are working in. 

Let $A$ be a dg algebra (or a dg category). We denote by $\mathcal D (A)
$ its (triangulated) derived category and by $\per (A)$ the \textit{perfect derived category}, that is, the subcategory of compact objects in $\mathcal D (A)$, which coincides with the smallest idempotent-complete triangulated subcategory of $\mathcal D (A)$ containing $A$. We denote by $\dba$ the \textit{bounded derived category} of $A$ which is the full subcategory of $\mathcal D (A)$ consisting of all dg $A$-modules with finite dimensional total cohomology. Recall that for a proper  dg algebra $A$ (that is, for $A \in \dba$), the fully faithful embedding $\per (A) \subseteq \dba$ is an equivalence  when $A$ is homologically smooth. In this case, if $A$ is zero-graded, then every finitely generated $A$ module has a finite projective resolution and can thus be modelled by a bounded complex of finitely generated projectives.

The triangulated subcategory of $\per(A)$ generated by the
free $A$-module $A$ has a natural dg enhancement given by the
category $\tw(A)$
of {\it twisted complexes} over the dg category  $\add A$, which is the additive closure of $A$ as a dg $A$-module \cite{BondalKapranov}. 
The category $\tw (\mathcal A)$ of twisted complexes over a dg  or $A_\infty$-category $\mathcal A$ can be viewed as the closure of $\mathcal A$ under shifts and mapping cones. This ensures that $\tw (\mathcal A)$ is a pre-triangulated dg  or $A_\infty$-category and its homotopy category $\mathrm H^0 \tw (\mathcal A)$ is triangulated and equivalent, up to direct factors, to $\per (\mathcal A)$. 

If $A$ is proper but not homologically smooth, the embedding $\per (A) \subseteq \dba$ is usually not essentially surjective. At least for the algebras of interest in this survey, we may use unbounded generalisations of twisted complexes in \cite{AL} to give an enhancement of $\dba$. Concretely, if $A$ is a finite dimensional graded gentle algebra (viewed as a proper dg algebra with trivial differential), its bounded derived category $\dba$ is generated by the simple $A$-modules \cite{BoothGoodbodyOpper}. 
In this case, a dg enhancement of $\dba$ can be given by the category generated by the (possibly unbounded) twisted complexes corresponding to the (possibly infinite) resolutions of the simples, viewed as a full subcategory of  $\Tw_\pm (\add A)$, the category of unbounded two-sided twisted complexes.

The theory of (possibly unbounded) twisted complexes thus provides a natural dg enhanced framework generalizing the classical equivalences $\per (A) \simeq \mathcal K^b (\proj A)$ and $\dba \simeq \mathcal K^{-,b} (\proj A)$ for a finite dimensional zero-graded algebra $A$, where $\mathcal K^b(\proj A)$ and $\mathcal K^{-,b}(\proj A)$ denote the homotopy categories of bounded complexes resp.~bounded above complexes of projectives with bounded cohomology.

\subsection{Geometric model for the bounded derived category of a graded gentle algebra}\label{sec: correspondence gentle algebra surface}

\begin{figure}[b]
    \centering
    \begin{tikzpicture}
        \begin{scope} % left
            \begin{scope}[xshift=-9em, yshift=4.75em, x=3em, y=1em] % algebra
            \node[font=\small, left=.3em, overlay] (A) at (0,0) {\strut$A \colon$};
            \node[font=\small, overlay] (1) at (0,0) {\strut$1$};
            \node[font=\small, overlay] (2) at (1,0) {\strut$2$};
            \node[font=\small, overlay] (3) at (2,0) {\strut$3$};
            \node[font=\small, ellipse, inner sep=0, outer sep=0pt, overlay] (4) at (3,0) {\strut$4\,,$};
            \path[-stealth, overlay] (1) edge node[above=-.2ex, font=\scriptsize, overlay] {$a$} (2) (2) edge node[above=-.2ex, font=\scriptsize, overlay] {$b$} (3) (3) edge[bend left=25, overlay] node[above=-.2ex, font=\scriptsize, overlay] {$c$} (4) (4) edge[bend left=25] node[below=-.2ex, font=\scriptsize] {$d$} (3);
            \node[font=\small, right=.5em, overlay] (I) at (3,0) {\strut$I = \langle dc, cd \rangle$};
            \draw[line width=.6pt, line cap=round, dash pattern=on 0pt off 1.2pt] (2,0) + (17:1.15em) arc[start angle=18, end angle=-19, radius=1.15em];
            \draw[line width=.6pt, line cap=round, dash pattern=on 0pt off 1.2pt, overlay] (3,0) + (180+17:1.15em) arc[start angle=180+18, end angle=180-19, radius=1.15em];
            \end{scope}
            \begin{scope}[xshift=-8em, yshift=-1em] % L
            \draw[fill=black] (0,0) circle(.1em);
            \draw[fill=black] (65:3em) circle(.1em);
            \draw[fill=black] (25:3em) circle(.1em);
            \draw[fill=black] (-25:3em) circle(.1em);
            \draw[line width=.5pt] (65:3em) to (0, 0) to (25:3em) (-25:3em) to[bend right=30] (0, 0) to[bend right=30] (-25:3em);
            \draw[-stealth] (63:1.1em) arc[start angle=63, end angle=27, radius=1.1em]; %to[bend left=20] (28:1.1em);
            \draw[-stealth] (23:1.1em) arc[start angle=23, end angle=-2, radius=1.1em]; %to[bend left=20] (-1:1.1em) to (-3:1.1em);
            \draw[-stealth] (-6:1.1em) arc[start angle=-6, end angle=-44, radius=1.1em]; %to[bend left=20] (-45:1.1em);
            \draw[-stealth] ($(-25:3em)+(180-6:1.1em)$) arc[start angle=180-6, end angle=180-44, radius=1.1em]; %$(-38:1.9em) to[bend left=30] (-13:1.9em);
            \node[font=\tiny] at (72:2.5em) {$1$};
            \node[font=\tiny] at (32:2.5em) {$2$};
            \node[font=\tiny] at (-10:2.5em) {$3$};
            \node[font=\tiny] at (-40:2.5em) {$4$};
            \node[font=\scriptsize] at (45:1.5em) {$a$};
            \node[font=\scriptsize] at (8:1.5em) {$b$};
            \node[font=\scriptsize] at (-23:1.45em) {$c$};
            \node[circle, inner sep=0, font=\scriptsize,fill=white] at (-45:1.7em) {$d$};
            \node[font=\tiny] at (65:3.5em) {$e_1$};
            \node[font=\tiny] at (25:3.55em) {$e_2$};
            \node[font=\tiny] at (-25:3.5em) {$d$};
            \node[font=\tiny] at (190:.75em) {$abc$};
            \end{scope}
            \begin{scope} % R
            \draw[line width=.5pt] (0,0) circle(3em);
            \draw[line width=.6pt, draw=yellow!30!green, line cap=round] (130:3em) to[bend left=25] (210:3em) to[bend right=8] (50:3em) (-30:3em) to[bend right=25] (210:3em) to[bend right=20] (-30:3em);
            \node[circle,fill=white,draw=black,inner sep=0,outer sep=0,minimum size=.35em] at (50:3em) {};
            \node[circle,fill=white,draw=black,inner sep=0,outer sep=0,minimum size=.35em] at (130:3em) {};
            \node[circle,fill=white,draw=black,inner sep=0,outer sep=0,minimum size=.35em] at (210:3em) {};
            \node[circle,fill=white,draw=black,inner sep=0,outer sep=0,minimum size=.35em] at (-30:3em) {};
            \node[circle,fill=green!10!red,inner sep=0,outer sep=0,minimum size=.4em] at (270:3em) {};
            \node[circle,fill=green!10!red,inner sep=0,outer sep=0,minimum size=.4em] at (170:3em) {};
            \node[circle,fill=green!10!red,inner sep=0,outer sep=0,minimum size=.4em] at (270:1.5em) {};
            \node[circle,fill=green!10!red,inner sep=0,outer sep=0,minimum size=.4em] at (10:3em) {};
            \node[circle,fill=green!10!red,inner sep=0,outer sep=0,minimum size=.4em] at (90:3em) {};
            \draw[line width=.6pt, draw=green!10!red] (270:3em) to (270:1.5em) to[bend left=15] (10:3em) to[bend left=20] (90:3em) to[bend left=20] (170:3em);
            \node[font=\footnotesize, color=green!10!red, right] at (3.3em,.5em) {$\Delta^*$};
            \node[font=\footnotesize, color=yellow!30!green, right] at (3.3em,-1.6em) {$\Delta$};
            \node[font=\tiny, color=yellow!30!green] at (180:2.3em) {$1$};
            \node[font=\tiny, color=yellow!30!green] at (90:1em) {$2$};
            \node[font=\tiny, color=yellow!30!green] at (-25:2.1em) {$3$};
            \node[font=\tiny, color=yellow!30!green] at (-61:2.5em) {$4$};
            \end{scope}
        \end{scope}
        \begin{scope}[xshift=24em] % right
            \begin{scope}[xshift=-10em, yshift=4.75em, x=3em, y=1em] % algebra
            \node[font=\small, left=.3em, overlay] (A) at (0,0) {\strut$A \colon$};
            \node[font=\small, overlay] (1) at (0,0) {\strut$1$};
            \node[font=\small, overlay] (2) at (1,0) {\strut$2$};
            \node[font=\small, overlay] (3) at (2,0) {\strut$3$};
            \node[font=\small, ellipse, inner sep=0, outer sep=0pt, overlay] (4) at (3,0) {\strut$4\,,$};
            \path[stealth-, overlay] (1) edge node[above=-.2ex, font=\scriptsize, pos=.6, overlay] {$a'$} (2) (2) edge node[above=-.2ex, font=\scriptsize, pos=.6, overlay] {$b'$} (3) (3) edge[bend left=25, overlay] node[above=-.2ex, font=\scriptsize, overlay] {$c'$} (4) (4) edge[bend left=25, overlay] node[below=-.2ex, font=\scriptsize] {$d'$} (3);
            \node[font=\small, right=.5em, overlay] (I) at (3,0) {\strut$I = \langle a' b', b' c' \rangle$};
            \draw[line width=.6pt, line cap=round, dash pattern=on 0pt off 1.2pt, overlay] (1,0) ++(-1em, -.1em) to[bend right=45, overlay] ++(2.2em,0) (2,0) +(-.8em, .1em) to[bend left=45, overlay] +(25:1.2em);
            \end{scope}
            \begin{scope}[xshift=-7.25em, yshift=.75em] % L
            \draw[fill=black] (180:2em) circle(.1em);
            \draw[fill=black] (90:2em) circle(.1em);
            \draw[fill=black] (0:2em) circle(.1em);
            \draw[fill=black] (270:2em) circle(.1em);
            \draw[fill=black] (270:3.5em) circle(.1em);
            \node[font=\tiny, overlay] at (270:4.1em) {$e_4$};
            \node[font=\tiny] at (135:1.8em) {$1$};
            \node[font=\tiny] at (45:1.8em) {$2$};
            \node[font=\tiny] at (-45:1.8em) {$3$};
            \node[font=\tiny, left=-.15em] at (270:3.1em) {$4$};
            \draw[stealth-] (90:2em) + (-130:.75em) arc[start angle=-130, end angle=-48, radius=.75em];
            \draw[stealth-] (0:2em) + (140:.75em) arc[start angle=140, end angle=222, radius=.75em];
            \draw[-stealth] (270:2em) + (265:.7em) arc[start angle=265, end angle=48, radius=.7em];
            \draw[-stealth] (270:2em) + (40:.7em) arc[start angle=40, end angle=-87, radius=.7em];
            \draw[line width=.5pt] (180:2em) to (90:2em) to (0:2em) to (270:2em) to (270:3.5em);
            \node[font=\tiny] at (180:2.6em) {$e_1$};
            \node[font=\tiny] at (90:2.6em) {$a'$};
            \node[font=\tiny] at (0:2.6em) {$b'$};
            \node[font=\tiny, ellipse, fill=white, inner sep=0, right=.15em] at (270:2em) {$cd$};
            \node[font=\scriptsize] at (90:.8em) {$a'$};
            \node[font=\scriptsize] at (-8:.87em) {$b'$};
            \node[font=\scriptsize, right=.1em] at (223:2.2em) {$c'$};
            \node[font=\scriptsize, right=.1em] at (270:2.8em) {$d'$};
            \end{scope}
            \begin{scope} % R
            \draw[line width=.5pt] (0,0) circle(3em);
            \draw[line width=.6pt, draw=yellow!30!green] (270:3em) to (270:1.5em) to[bend left=15] (10:3em) to[bend left=20] (90:3em) to[bend left=20] (170:3em);
            \draw[line width=.6pt, draw=green!10!red, line cap=round] (130:3em) to[bend left=25] (210:3em) to[bend right=8] (50:3em) (-30:3em) to[bend right=25] (210:3em) to[bend right=20] (-30:3em);
            \node[circle,fill=white,draw=black,inner sep=0,outer sep=0,minimum size=.35em] at (270:3em) {};
            \node[circle,fill=white,draw=black,inner sep=0,outer sep=0,minimum size=.35em] at (170:3em) {};
            \node[circle,fill=white,draw=black,inner sep=0,outer sep=0,minimum size=.35em] at (270:1.5em) {};
            \node[circle,fill=white,draw=black,inner sep=0,outer sep=0,minimum size=.35em] at (10:3em) {};
            \node[circle,fill=white,draw=black,inner sep=0,outer sep=0,minimum size=.35em] at (90:3em) {};
            \node[circle,fill=green!10!red,inner sep=0,outer sep=0,minimum size=.4em] at (50:3em) {};
            \node[circle,fill=green!10!red,inner sep=0,outer sep=0,minimum size=.4em] at (130:3em) {};
            \node[circle,fill=green!10!red,inner sep=0,outer sep=0,minimum size=.4em] at (210:3em) {};
            \node[circle,fill=green!10!red,inner sep=0,outer sep=0,minimum size=.4em] at (-30:3em) {};
            \node[font=\footnotesize, color=yellow!30!green, right] at (3.3em,.5em) {$\Delta$};
            \node[font=\footnotesize, color=green!10!red, right] at (3.3em,-1.6em) {$\Delta^*$};
            \node[font=\tiny, color=yellow!30!green, left=-.2em] at (270:2.5em) {$4$};
            \node[font=\tiny, color=yellow!30!green, left=-.2em] at (-10:2.2em) {$3$};
            \node[font=\tiny, color=yellow!30!green] at (70:2.5em) {$2$};
            \node[font=\tiny, color=yellow!30!green] at (110:2.5em) {$1$};
            \end{scope}
        \end{scope}
    \end{tikzpicture}
    \caption{Two mutually Koszul dual (graded) gentle algebras, their ribbon graphs and the corresponding admissible dissections. Every path in the quiver induces a morphism between the corresponding indecomposable  dg $A$-modules $e_i A = P_i$. For example, for the model on the left side, we have morphisms $P_1 \stackrel{a.}\longrightarrow P_2[|a|]$, $P_2 \stackrel{b.}\longrightarrow P_3[|b|]$ and $P_1 \stackrel{ba.}\longrightarrow P_3[|a|+|b|]$, where $|a|$ and $|b|$ denote the degrees of $a$ and $b$, respectively.}
     \label{fig:dual dissections}
\end{figure}
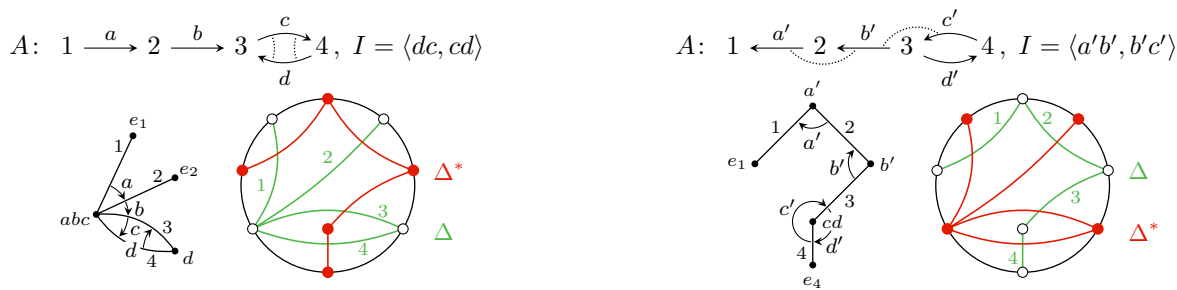

We now introduce the surface model $S$ for the bounded derived category $\dba$ of a graded gentle algebra $A$ given in \cite{OpperPlamondonSchroll}. As described in the last paragraph of Section \ref{sec: geometric model intro}, the surface $S$ arises from the ribbon graph $\Delta$ associated to $A$ as in \cite{Schroll15}. Concretely, given a finite dimensional (graded) gentle algebra $A = \Bbbk Q / I$, the vertices of its associated ribbon graph $\Delta$ are given by the maximal paths in $(Q,I)$, as well as all length zero paths whose vertices appear only once in all maximal paths. Each vertex $v$ of $Q$ gives rise to an edge of $\Delta$, joining the two vertices of $\Delta$ corresponding to the paths in which $v$ appears. The arrows  in a  path corresponding to a vertex of $\Delta$ give rise to a linear order, whose cyclic completion gives $\Delta$ the structure of a ribbon graph. The construction of $\zD$ for infinite dimensional graded gentle algebras is adapted accordingly, see Figure~\ref{fig:dual dissections} for an example.

A \textit{surface dissection} is given by a collection of curves in $S$ cutting the surface into polygons. Here we consider the surface dissection induced by the embedding of the ribbon graph $\zD$ into $S$, where the vertices of $\Delta$ are glued to the  boundary of $S$. The vertices correspond to marked points denoted $\gpoint$ and the curves corresponding to the edges of  $\Delta$ are depicted in green in the figures. The dual of $\Delta$, denoted $\Delta^*$ (depicted in red), then serves as a kind of coordinate system  for $\dba$. Indecomposable objects of $\dba$ correspond to certain curves in $S$ whose representation as twisted complexes can be read directly from the crossings with curves belonging to $\Delta^*$.
The vertices $\gpoint$ of $\zD$ and the vertices $\rpoint$ of $\zD^*$ alternate on every boundary component.  For homologically non-smooth gentle algebras, there are, in addition, $\rpoint$-points in the interior of the surface. As a consequence, the dissection  $\zD$ is such  that each polygon contains exactly one $\rpoint$-point, either in its boundary or in its interior. We call any dissection with this property an \textit{admissible dissection}. See Figure \ref{fig:dual dissections} for  concrete examples. 

Furthermore, the grading of the gentle algebra induces a grading of the ribbon graph $\zD$, this in turn induces a grading on all curves in the associated surface $S$. It turns out that for closed curves,  this grading corresponds exactly to the winding number of the closed curves. Here the winding numbers are those with respect to  a line field on the surface constructed from the graded gentle algebra in \cite{LekiliPolishchuk}.

\paragraph{A note on Koszul duality} Koszul duality of graded gentle algebras swaps the roles of $\zD$ and $\zD^*$. In particular, we note that a proper homologically non-smooth graded gentle algebra has $\rpoint$-points in the interior of the surface  (corresponding to cycles with relations at every subpath of length 2) but no $\gpoint$-points. In contrast, a homologically smooth but non-proper graded gentle algebra will have $\gpoint$-points in the interior of the surface (corresponding to cycles with no relations) but no interior $\rpoint$-points. See, for example, the two Koszul dual algebras in Figure~\ref{fig:dual dissections}. The surface of a  non-smooth and non-proper graded gentle algebra will have both $\rpoint$-points and $\gpoint$-points in the interior of the surface whose roles are exchanged by Koszul duality.

\subsection{Bijection between indecomposable objects in $\dba$ and graded curves on the surface} In \cite{HaidenKatzarkovKontsevich} it is shown that for a  graded gentle algebra $A$ there is a bijection between the indecomposable objects (in terms of twisted complexes) in $\ka$ and homotopy classes of curves with local system in the associated graded marked surface. Coming from a representation theoretic perspective, the same has been shown in \cite{OpperPlamondonSchroll} for $\dba$. For the  zero-graded case, this builds on the classification of indecomposable objects in terms of homotopy strings and bands \cite{BekkertMarcosMerklen}, see also \cite{BurbanDrozd}. In the general graded case it additionally uses \cite{HaidenKatzarkovKontsevich} and \cite{AL}  as well as Koszul duality. We note that the new graded version of \cite{OpperPlamondonSchroll}  also relies on the proof in \cite{BoothGoodbodyOpper}  that, for graded gentle algebras, the thick subcategory of the derived category generated by the simple objects coincides with $\dba$.

 \begin{figure}
    \centering
     \begin{tikzpicture}
        \begin{scope}[xshift=-2em, yshift=4.75em, x=3em, y=1em] % algebra
        \node[font=\small, left=.3em, overlay] (A) at (0,0) {\strut$A \colon$};
        \node[font=\small, overlay] (1) at (0,0) {\strut$1$};
        \node[font=\small, overlay] (2) at (1,0) {\strut$2$};
        \node[font=\small, overlay] (3) at (2,0) {\strut$3$};
        \node[font=\small, ellipse, inner sep=0, outer sep=0pt, overlay] (4) at (3,0) {\strut$4\,,$};
        \path[-stealth, overlay] (1) edge node[above=-.2ex, font=\scriptsize, overlay] {$a$} (2) (2) edge node[above=-.2ex, font=\scriptsize, overlay] {$b$} (3) (3) edge[bend left=25, overlay] node[above=-.2ex, font=\scriptsize, overlay] {$c$} (4) (4) edge[bend left=25] node[below=-.2ex, font=\scriptsize] {$d$} (3);
        \node[font=\small, right=.5em, overlay] (I) at (3,0) {\strut$I = \langle dc, cd \rangle$};
        \draw[line width=.6pt, line cap=round, dash pattern=on 0pt off 1.2pt] (2,0) + (17:1.15em) arc[start angle=18, end angle=-19, radius=1.15em];
        \draw[line width=.6pt, line cap=round, dash pattern=on 0pt off 1.2pt, overlay] (3,0) + (180+17:1.15em) arc[start angle=180+18, end angle=180-19, radius=1.15em];
        \end{scope}
        \begin{scope}[xshift=2.5em] % right surface
        \draw[line width=.5pt, draw=yellow!30!green, line cap=round, opacity=.2] (130:3em) to[bend left=25] (210:3em) to[bend right=8] (50:3em) (-30:3em) to[bend right=25] (210:3em) to[bend right=20] (-30:3em);
        \node[font=\tiny, color=green!10!red, opacity=1] at (110:2.5em) {$1$};
        \node[font=\tiny, color=green!10!red, opacity=1] at (70:2.5em) {$2$};
        \node[font=\tiny, color=green!10!red, opacity=1, left=-.2em] at (-10:2.2em) {$3$};
        \node[font=\tiny, color=green!10!red] at (263:2.55em) {$4$};
        \draw[line width=.5pt] (0,0) circle(3em);
        \draw[line width=.6pt, draw=green!10!red, opacity=1] (270:3em) to (270:1.5em) to[bend left=15] (10:3em) to[bend left=20] (90:3em) to[bend left=20] (170:3em);
        \draw[line width=.6pt, draw=blue, line cap=round] (130:3em) to[bend right=8] (-30:3em) (130:3em) to[in=180, out=-60] ($(270:1.5em)+(90:.85em)$) arc[start angle=90, end angle=-135, radius=.8em] arc[start angle=225, end angle=45, radius=.7em] arc[start angle=45, end angle=-135, radius=.6em] arc[start angle=225, end angle=45, radius=.5em] arc[start angle=45, end angle=-135, radius=.4em];
        \draw[stealth-] (158:1.3em) arc[start angle=-70, end angle=-45, radius=1.6em];
        \node[font=\scriptsize] at (158:.75em) {$f$};
        \node[font=\scriptsize, color=blue] at (106:1.2em) {$\alpha$};
        \node[font=\scriptsize, color=blue] at (169:1.65em) {$\beta$};
        \draw[line width=.7pt, draw=blue, line cap=round, dash pattern=on 0pt off 1.2pt] ($(270:1.5em)+(-155:.35em)$) arc[start angle=-155, end angle=-220, radius=.3em];
        \node[circle,fill=white,draw=black,inner sep=0,outer sep=0,minimum size=.35em] at (50:3em) {};
        \node[circle,fill=white,draw=black,inner sep=0,outer sep=0,minimum size=.35em] at (130:3em) {};
        \node[circle,fill=white,draw=black,inner sep=0,outer sep=0,minimum size=.35em] at (210:3em) {};
        \node[circle,fill=white,draw=black,inner sep=0,outer sep=0,minimum size=.35em] at (-30:3em) {};
        \node[circle,fill=green!10!red,inner sep=0,outer sep=0,minimum size=.4em] at (270:3em) {};
        \node[circle,fill=green!10!red,inner sep=0,outer sep=0,minimum size=.4em] at (170:3em) {};
        \node[circle,fill=green!10!red,inner sep=0,outer sep=0,minimum size=.4em] at (270:1.5em) {};
        \node[circle,fill=green!10!red,inner sep=0,outer sep=0,minimum size=.4em] at (10:3em) {};
        \node[circle,fill=green!10!red,inner sep=0,outer sep=0,minimum size=.4em] at (90:3em) {};
        \end{scope}
        \begin{scope}[xshift=6em, yshift=2em]
        \node[font=\footnotesize, align=left, anchor=north, right, text depth=3.1em] at (17.9em, -3em) {$d_\beta = $};
        \node[font=\footnotesize, align=left, anchor=north, right] at (1em,0) {$X_\alpha = P_1 \oplus P_3 [|a| + |b| - 1]$};
        %\node[font=\footnotesize, align=left, anchor=north, right] at (0,-4em) {$X_\delta = P_1 \oplus P_3 [|a| + |b| - 1]$ \\ $\qquad\qquad \oplus\bigoplus\limits_{k \geq 1} \Bigl( \begin{pmatrix} P_4 [|a| + |b| - 1 - (k-1) (|c|-1) - k (|d| - 1)] \\ \oplus P_3 [|a| + |b| - 1 - k (|c|-1) - k (|d| - 1)] \end{pmatrix} \Bigr)$};
        \node[font=\footnotesize, align=left, anchor=north, right] at (1em,-3em) {$X_\beta = P_1 \oplus P_3 [|a| + |b| - 1]$ \\ $\phantom{X_\beta = P_1} {}\oplus P_4 [|a| + |b| - 1 - |d| + 1]$ \\ $\phantom{X_\beta = P_1} {}\oplus P_3 [|a| + |b| - 1 - |d| + 1 - |c| + 1]$ \\ $\phantom{X_\beta = P_1} {}\oplus \dotsb$};
        \begin{scope}[xshift=-1.9em, yshift=.3em]
        \node[right, font=\large] at (21.7em,-3em) {$\Biggl( \hspace{3.6em} \Biggr)$};
        \node[matrix of math nodes, every node/.style={font=\scriptsize}, ampersand replacement=\&, row sep={.8em,between origins}, column sep={.8em, between origins}] at (25em, -3em) {$0$ \& $0$ \& $0$ \& $0$ \& ${\cdot}{\cdot}{\cdot}$ \\ $ba$ \& $0$ \& $d$ \& $0$ \& \\ $0$ \& $0$ \& $0$ \& $c$ \& \\ $\raisebox{2.1pt}{$\cdot$}\mkern-5.2mu{\cdot}\mkern-5.2mu\raisebox{-2.1pt}{$\cdot$}$ \& \& \& \& $\raisebox{1.5pt}{$\cdot$}\mkern-1mu{\cdot}\mkern-1mu\raisebox{-1.5pt}{$\cdot$}$ \\};
        \node[right, font=\large] at (22.08em,0) {$\bigl( \hspace{1.4em} \bigr)$};
        \node[left, font=\footnotesize] at (22.58em,-.31em) {$d_\alpha = {}$};
        \node[matrix of math nodes, every node/.style={font=\scriptsize}, ampersand replacement=\&, row sep={.8em,between origins}, column sep={.8em, between origins}] at (23.8em, 0em) {$0$ \& $0$ \\ $ba$ \& $0$ \\};
        \end{scope}
        \end{scope}
     \end{tikzpicture}
    %\includegraphics[width=0.5\linewidth]{}
   % \caption{Caption}
    \caption{Example of a bounded twisted complex $X_\alpha$ and an unbounded twisted complex $X_\beta$ with differentials $d_\alpha$ and $d_\beta$, respectively, and a morphism $f$ between them induced by the identity on the component $P_1 \oplus P_3 [|a|+|b|-1]$.}
     \label{fig: Example Twisted Complex}
\end{figure}

\paragraph{Bijection between indecomposable objects and curves} 
 The general idea is that the edges of $\zD^*$ serve  as a kind of coordinate system. In its simplest form, the curves given by the edges of $\zD$  correspond to indecomposable  dg $A$-modules   of the form $P_i=e_iA$, where $i$ is the vertex of $Q$ corresponding to the edge of $\zD$,  
see Figure~\ref{fig:dual dissections}.
More generally, given a curve $\gamma$, the corresponding twisted complex $X_\gamma$ with differential $d_\gamma$ is  constructed from the successive crossings of $\gamma$ with the edges of $\zD^*$. The differential $d_\gamma$ is induced by the corresponding paths in $Q$. We illustrate this  in Figure~\ref{fig: Example Twisted Complex}.
We note that, in general, the gradings on curves, which also arise from the crossings of the curve with the dual dissection $\zD^*$, encode the action of the  shift functor on the objects in $\dba$.  

To summarise,  the twisted complexes associated to  homotopy strings correspond on the surface to  homotopy classes of open curves (connecting $\gpoint$-points) equipped with a grading.  The twisted complexes associated to homotopy bands correspond to  homotopy classes of (gradable) closed curves, see  \cite{OpperPlamondonSchroll} for details. In terms of the line field defined from a graded gentle algebra in  \cite{LekiliPolishchuk}, a gradable closed curve is a curve with winding number zero.  We will refer to the homotopy classes of curves connecting $\gpoint$-points as \textit{graded arcs} and to the homotopy classes of gradable closed curves as \textit{graded closed curves}.  Just as homotopy bands come together with a finite dimensional indecomposable $\Bbbk [X^\pm]$-module, the associated graded closed curves come with an indecomposable local system of finite rank. This describes the indecomposable objects in $\ka$.

To pass from $\ka$ to $\dba$, we need to add additional curves. For a non-smooth graded gentle algebra, the associated surface will have $\rpoint$-points in its interior. In order to establish a correspondence between indecomposable objects in $\dba$ and curves on the surface, we need to add to the set of graded arcs and curves described above, all (homotopy classes of) graded curves which spiral at either or both ends (locally in the clockwise direction) around a $\rpoint$-point on the surface together with a grading.

\subsection{Geometric interpretation of morphisms and mapping cones in $\dba$ } A basis of the morphism space between two indecomposable objects in $\dba$, for $A$ a finite dimensional gentle algebra, is given explicitly in \cite{ArnesenLakingPauksztello} using (homotopy) string combinatorics. It is shown in \cite{OpperPlamondonSchroll} that 
the results in \cite{ArnesenLakingPauksztello}
still hold in the graded case and thus give an algebraic description of a basis of the morphism space between indecomposable objects in $\dba$, for $A$ a proper graded gentle algebra.
Based on the description of these basis morphisms in \cite{ArnesenLakingPauksztello} and their generalisation to the graded case, it is shown in \cite{OpperPlamondonSchroll} that there is a bijection between basis morphisms and (oriented) intersections in the graded marked surface associated to $A$.

Furthermore, the mapping cones of the basis morphisms  can be read from the surface model. Namely, given a basis morphism corresponding to an oriented intersection, its mapping cone is given by the direct sum of the indecomposables associated to the curves resolving the crossing.  This result is based on the calculation of these mapping cones in terms of homotopy string combinatorics in \cite{CanakciPauksztelloSchrollMappingCone} which again can be generalised to the graded case.  

%, see Figure~\ref{Fig:Mapping Cone}. 
%\newline
\begin{wrapfigure}{r}{0.18\textwidth}
\centering 
\begin{tikzpicture}%[transform canvas={scale=0.7}]
[>=stealth,x=1.45em, y=1.45em]
%\draw (-0.5,0) node[black!60] {$\aaa$};
%\draw [bend left,line width=.5pt,->] (-0.3,0.3) to (0.3,0.3); 
\draw (0,0.97) node[font=\footnotesize] {$\aaa$};
\draw[line width=.5pt,->, overlay] (0.4,-0.45) arc[start angle=-45, end angle=-135, radius=.56];
\draw[line width=.5pt,<-, overlay] (0.4,0.45) arc[start angle=45, end angle=135, radius=.56];
\draw (0,-0.97) node[font=\footnotesize] {$\aaa$};
\draw[black, line width=.5pt, ->] (-2,2)  --  (2,-2);
\draw (-1.5,1) node[font=\footnotesize] {$\za$};
\draw[black, line width=.5pt, <-] (-2,-2) -- (2,2);
\draw (1.5,1) node[font=\footnotesize] {$\zb$};
\draw[line width=.5pt, ->, overlay] (-1.8,2.2) to[bend right=40, overlay] (1.8,2.2);
\draw (0,1.9) node[overlay, font=\footnotesize] {$\gamma_1$};
\draw[line width=.5pt, <-, overlay] (-1.8,-2.2) to[bend left=40, overlay] (1.8,-2.2);
\draw (0,-1.9) node[overlay, font=\footnotesize] {$\gamma_2$};
%\filldraw[black] (0,0) circle (2pt) node[anchor=west]{Intersection point};
\end{tikzpicture}
%%\caption{}
%\label{Fig:Mapping Cone}
\end{wrapfigure}  
More precisely, for $X_\alpha$ and $X_\beta$ twisted complexes associated to the curves $\za$ and $\zb$, let $\aaa: X_\za \to X_\zb$ be the morphism associated to the oriented intersection $\aaa$. Then there exist gradings on $\zg_1$ and $\zg_2$ such that there is a triangle $X_\za \smash{\stackrel{\aaa}\to} X_\zb \to X_{\zg_1} \oplus X_{\zg_2}$, where $X_{\zg_1}$ and $X_{\zg_2}$ are the twisted complexes associated to the graded curves $\zg_1$ and $\zg_2$.
We note that either of $X_{\zg_1}$ or $X_{\zg_2}$ might be zero if the intersection of $\za$ and $\zb$ lies in the boundary, or that $\zg_1$ might be identified with $\zg_2$ if the resolution of the crossing corresponds to a closed curve. In the latter case the corresponding object $X_{\gamma_1}$ might be decomposable.

\subsection{Geometric interpretation of the Serre functor in $\ka$}\label{sec: Serre functor} Since a proper graded gentle algebra $A$ is Gorenstein, the perfect derived  category $\ka$ has a Serre functor $\SSS$ given by the derived Nakayama functor $\nu = - \otimes^{\mathbb{L}}_A DA$, see, for example \cite{ Jin, KellerDerivingDG, KrauseBook}. This is the shift of the  Auslander--Reiten translate $\tau$, that is $\SSS =  \tau[1]$. Furthermore, it is shown in \cite{OpperPlamondonSchroll}, that given a graded arc $\zg$ with associated twisted complex $X_\zg$ connecting $\gpoint$-points on $S$, then $\SSS (X_\zg)$ is the twisted complex corresponding the curve $\zg'$ obtained from $\zg$ by moving each of the endpoints of the curve to the next $\gpoint$-point  on the same boundary component in the clockwise direction.  On twisted complexes corresponding to graded closed curves, the Serre functor acts as the  shift only changing the grading of the curve.  

The Serre functor in derived categories and its categorical entropy, defined in \cite{DHKK}, have been intensively studied, see for example, \cite{Elagin, ElaginLunts}. Recently, in \cite{ChangElaginSchroll} the categorical entropy of the Serre functor in the perfect derived category of a proper, homologically smooth graded gentle algebra has  been given in terms of  winding numbers of boundary components and their number of marked points in the associated graded marked surface.

\subsection{Recollements and surface cuts}

In this section we see that recollements of bounded derived categories of graded gentle algebras can be expressed in terms of surface cuts, see also  \cite{Dyckerhoff}. This provides reduction techniques which are used in the proofs of many of the results presented below. We note that a similar idea of cutting  surfaces was  used in the context of Bridgeland stability conditions for Fukaya categories of surfaces in \cite{Takeda}.

In \cite{ChangJinSchroll}, given a homologically smooth and proper graded gentle algebra $A(\zD)$ with associated surface dissection $\zD$, it is shown that cutting (and contracting along the cut)  a set $\zG$ of arcs of $\zD$, gives rise to a recollement of bounded derived categories of graded gentle algebras. The cut set $\zG$ is itself a ribbon graph (possibly not connected) with  grading induced by the grading of $\zD$ and gives rise to a gentle algebra $A(\zG)$.
On the other hand, the triangulated category of the left side of the recollement is equivalent to the bounded derived category of a graded gentle algebra $A_{\zG}$. 
Then by \cite{ChangJinSchroll} there is a recollement of triangulated categories (see  Figure~\ref{fig:choice of diss})
$$
\begin{tikzpicture}
   % \begin{scope}[shift={(5.0,0)}]
  \node (c1) at (0,0) {$\mathcal{D}^b(A_\zG)$};
  \node (c2) at (3,0) {$\mathcal{D}^b(A(\zD))$};
  \node (c3) at (6,0) {$\mathcal{D}^b(A(\zG)).$};
  \draw[stealth-, overlay, transform canvas={yshift=.9ex}] (c1.east) to[bend left=20] (c2.west);
  \draw[-stealth] (c1) to (c2);
  \draw[stealth-, overlay, transform canvas={yshift=-.9ex}] (c1.east) to[bend right=20] (c2.west);
  \draw[stealth-, overlay, transform canvas={yshift=.9ex}] (c2.east) to[bend left=20] (c3.west);
  \draw[-stealth] (c2) to (c3);
  \draw[stealth-, overlay, transform canvas={yshift=-.9ex}] (c2.east) to[bend right=20] (c3.west);
\end{tikzpicture}
  $$
 The proof is based on a partial cofibrant dg algebra resolution of a graded quadratic monomial algebra relative to a subset of the quadratic relations.
A similar result for localising at a band object is considered in \cite{Bodin}.

%$$
%\begin{tikzpicture}
%   % \begin{scope}[shift={(5.0,0)}]
%  \node (c1) at (0,0) {1};
%  \node (c2) at (1.4,0) {2};
%  \node (c3) at (2.8,0) {3};
%  \node (c4) at (4.2,0) {4};
%  \node (c5) at (5.6,0) {5};
%  \node (c6) at (7,0) {6};
%  \node (c7) at (8.4,0) {7};
%  \node (c8) at (9.8,0) {8};
%  \node (c9) at (11.2,0) {9};
  
  %  labelled arrows 1 -> 9
%  \draw[->,overlay] (c1) to node[above,midway, overlay] {\small{$a_1$}}  (c2)
  %;
  %\draw[->,overlay] (c2) to node[above,midway, overlay] {\small{$a_2$}} (c3);
  %\draw[->,overlay] (c3) to node[above,midway, overlay] {\small{$a_3$}} (c4);
  %\draw[->,overlay] (c4) to node[above,midway, overlay] {\small{$a_4$}} (c5);
  %\draw[->,overlay] (c5) to node[above,midway, overlay] {\small{$a_5$}} (c6);
 % \draw[->,overlay] (c6) to node[above,midway, overlay] {\small{$a_6$}} (c7);
%  \draw[->,overlay] (c7) to node[above,midway, overlay] {\small{$a_7$}} (c8);
%  \draw[->,overlay] (c8) to node[above,midway, overlay] {\small{$a_8$}} (c9);
  
%  relations
%  \draw[densely dotted,overlay] (1,0.2) to [bend left] (1.8,0.2);
%  \draw[densely dotted,overlay] (2.4,0.2) to [bend left] (3.2,0.2);
%  \draw[densely dotted,overlay] (3.8,0.2) to [bend left] (4.6,0.2);

%\end{scope}

%\end{tikzpicture}
%$$
%where we cut along the edges $\zG=\{2,3,4,8\}$  of the associated admissible surface dissection $\zD$.

 \begin{figure}
\centering
\begin{tikzpicture}[x=1em, y=1em]
\begin{scope}[xshift=-11.5em] % LHS
\node[font=\small] at (0,-4.25em) {$\mathcal D^b (A_\Gamma)$};
\begin{scope}[xshift=-3em] % left circle
\draw[line width=.5pt] (0,0) circle(2em);
\draw[line width=.6pt, draw=yellow!30!green] (180:2em) to[bend left=25] (180+72:2em) to[bend left=25] (180+144:2em) (180-72:2em) to[bend left=5] (180+72:2em) to[bend left=5] (36:2em);
\foreach \a in {0,72,...,288} {
\node[circle,fill=green!10!red,inner sep=0,outer sep=0,minimum size=.4em] at (\a:2em) {};
\node[circle,fill=white,draw=black,inner sep=0,outer sep=0,minimum size=.35em] at (180+\a:2em) {};
}
\node[font=\tiny, color=yellow!30!green] at (195:1.1em) {$1$};
\node[font=\tiny, color=yellow!30!green] at (134:1.3em) {$5$};
\node[font=\tiny, color=yellow!30!green] at (46:1.3em) {$6$};
\node[font=\tiny, color=yellow!30!green] at (-38:1.35em) {$7$};
\end{scope}
\node at (0,0) {$\times$};
\begin{scope}[xshift=2em] % right circle
\draw[line width=.5pt] (0,0) circle(1em);
\draw[line width=.6pt, draw=yellow!30!green] (90:1em) to (270:1em);
\foreach \a in {0,180} {
\node[circle,fill=green!10!red,inner sep=0,outer sep=0,minimum size=.4em] at (\a:1em) {};
\node[circle,fill=white,draw=black,inner sep=0,outer sep=0,minimum size=.35em] at (90+\a:1em) {};
}
\node[font=\tiny, color=yellow!30!green] at (180:.3em) {$9$};
\end{scope}
\end{scope}
\begin{scope} % middle
 %\node[font=\small] at (0,4.75em) {$A (\Delta) : 1 \overset{a_1}\longrightarrow 2 \overset{a_2}\longrightarrow 3 \overset{a_3}\longrightarrow 4 \overset{a_4}\longrightarrow 5 \overset{a_5}\longrightarrow 6 \overset{a_6}\longrightarrow 7 \overset{a_7}\longrightarrow 8 \overset{a_8}\longrightarrow 9  $};
\node[font=\small] at (0,-4.25em) {$\mathcal D^b (A (\Delta))$};
\draw[line width=.6pt, draw=yellow!30!green] (180-54:3em) to[bend left=35] (180-18:3em) to[bend left=35] (180+18:3em) to[bend left=35] (180+54:3em) to[bend left=35] (180+90:3em) to (90:3em) (54:3em) to[bend right=5] (270:3em) to[bend left=10] (18:3em) (-18:3em) to[bend right=20] (270:3em) to[bend left=35] (270+36:3em);
\draw[line width=.5pt] (0,0) circle(3em);
\foreach \a in {0,36,...,324} {
\node[circle,fill=green!10!red,inner sep=0,outer sep=0,minimum size=.4em] at (\a:3em) {};
\node[circle,fill=white,draw=black,inner sep=0,outer sep=0,minimum size=.35em] at (\a+18:3em) {};
}
\node[font=\tiny, color=yellow!30!green] at (180-36:2.25em) {$1$};
\node[font=\tiny, color=yellow!30!green] at (180:2.25em) {$2$};
\node[font=\tiny, color=yellow!30!green] at (180+36:2.25em) {$3$};
\node[font=\tiny, color=yellow!30!green] at (180+72:2.25em) {$4$};
\node[font=\tiny, color=yellow!30!green] at (99:2.1em) {$5$};
\node[font=\tiny, color=yellow!30!green] at (58:2.2em) {$6$};
\node[font=\tiny, color=yellow!30!green] at (18:2.3em) {$7$};
\node[font=\tiny, color=yellow!30!green] at (-18:2.4em) {$8$};
\node[font=\tiny, color=yellow!30!green] at (304:2.5em) {$9$};
\draw[stealth-, transform canvas={yshift=.9ex}] (-7.5em,0) to[bend left=20] (-4em,0);
\draw[-stealth] (-7.5em,0) to (-4em,0);
\draw[stealth-, transform canvas={yshift=-.9ex}] (-7.5em,0) to[bend right=20] (-4em,0);
\draw[stealth-, transform canvas={yshift=.9ex}] (4em,0) to[bend left=20] (7.5em,0);a
\draw[-stealth] (4em,0) to (7.5em,0);
\draw[stealth-, transform canvas={yshift=-.9ex}] (4em,0) to[bend right=20] (7.5em,0);
\end{scope}
\begin{scope}[xshift=10.5em] % RHS
\node[font=\small] at (0,-4.25em) {$\mathcal D^b (A (\Gamma))$};
\node[font=\footnotesize, right] at (3em,0) {$\Gamma = \{ 2, 3, 4, 8 \}$};
\draw[line width=.6pt, draw=yellow!30!green] (108:2em) to[bend left=30] (180:2em) to[bend left=30] (252:2em) to[bend left=30] (-36:2em) to[bend left=30] (36:2em);
\draw[line width=.5pt] (0,0) circle(2em);
\foreach \a in {0,72,...,288} {
\node[circle,fill=green!10!red,inner sep=0,outer sep=0,minimum size=.4em] at (\a:2em) {};
\node[circle,fill=white,draw=black,inner sep=0,outer sep=0,minimum size=.35em] at (\a+36:2em) {};
}
\node[font=\tiny, color=yellow!30!green] at (144:.95em) {$2$};
\node[font=\tiny, color=yellow!30!green] at (216:.95em) {$3$};
\node[font=\tiny, color=yellow!30!green] at (288:.9em) {$4$};
\node[font=\tiny, color=yellow!30!green] at (20:1em) {$8$};
\end{scope}
\end{tikzpicture}
     \caption{Example of a recollement of  derived categories of graded gentle algebras in terms of surface cuts where the algebra $A(\zD)$ is given by
      $ 1 \tikztoarg{a_1} \rel 2 \tikztoarg{a_2} \rel 3 \tikztoarg{a_3} \rel 4
     \tikztoarg{a_4} 5 
     \tikztoarg{a_5} 6
     \tikztoarg{a_6} 7
     \tikztoarg{a_7} 8
     \tikztoarg{a_8} 9$ and where  $\Gamma = \{2,3,4,8\}$. Then $A_\zG \colon (1 \tikztoarg{b} 5 \tikztoarg{a_5} 6 \tikztoarg{a_6} 7) \times 9$ with $|b| = |a_1|+|a_2|+|a_3|+|a_4|-3$ and $A(\zG) \colon 2 \tikztoarg{a_2} \rel 3 \tikztoarg{a_3} \rel 4 \tikztoarg{c} 8$ with $|c| = |a_4|+|a_5|+|a_6|+|a_7|$.}
 \label{fig:choice of diss}
 \end{figure}

\section{Fukaya categories of surfaces with stops}\label{sec: Fukaya}

In symplectic geometry, categories of twisted complexes (see Section \ref{sec: twisted complexes}) are a central tool in the study of Fukaya categories \cite{SeidelBook}. In this context, one usually considers the category $\tw (\mathcal A)$ of bounded twisted complexes over an $A_\infty$-category $\mathcal A$. The objects of $\mathcal A$ correspond to a set of generating Lagrangians, with morphisms reflecting the Lagrangian intersection theory and higher morphisms corresponding to the count of pseudoholomorphic disks. Then  $\tw (\mathcal A)$ is  a pre-triangulated $A_\infty$-category, sometimes referred to as \textit{topological Fukaya category}. The category $\tw (\mathcal A)$ contains not only the generating Lagrangians in $\mathcal A$, but also twisted complexes such as mapping cones of morphisms which model connected sums (Lagrangian surgery) of Lagrangians \cite{FukayaOhOhtaOnoIII, FukayaOhOhtaOnoII}. 

Gentle algebras have a long history in representation theory (see Section \ref{sec: gentle}) and graded gentle algebras appear in \cite{HaidenKatzarkovKontsevich} as formal generators of topological Fukaya categories of marked surfaces. The surface model of Section \ref{sec: OPS Model} can be viewed as an exact symplectic surface and the curves belonging to the dissection $\Delta$ as a generating set of Lagrangian submanifolds. The $\rpoint$-points correspond to stops in the sense of partially wrapped Floer theory \cite{AurouxICM, Auroux,GanatraPardonShendeI, GanatraPardonShendeII}.
 Furthermore,   $\gpoint$-points and $\rpoint$-points in the interior of the surface  translated to the model of the topological Fukaya category in \cite{HaidenKatzarkovKontsevich},  correspond to unstopped and fully stopped boundary components, respectively.

This symplectic viewpoint sheds new light on the study of gentle algebras and also gives a strong motivation for considering gentle algebras with a $\mathbb Z$-grading and possibly a nontrivial $A_\infty$-structure. In this graded context, one can also consider polygons without boundary segments and no interior $\rpoint$-points, in addition to the polygons with $\rpoint$-points considered in Section \ref{sec: OPS Model}.  These new types of polygons can be understood as pseudoholomorphic disks, giving rise to a higher product $\mu_n$, equipping the corresponding gentle algebra with a nontrivial $A_\infty$-structure.
These polygons naturally arise when adding extra arcs to an admissible dissection as in Section~\ref{sec: OPS Model}. Namely,  for each $\rpoint$-point in the boundary one adds an arc, such that the $\rpoint$-point lies in a 2-gon. In the language of \cite{HaidenKatzarkovKontsevich} the resulting surface dissection is  an example 
of a \textit{full arc system}, see the left side of Figure~\ref{fig: HKK disk}.

Viewing the graded gentle algebra, equipped with this $A_\infty$-structure from polygons without $\rpoint$-points, as an $A_\infty$-category $\mathcal A$, one can define the (topological) partially wrapped Fukaya category of the graded marked surface $(S, M, \eta)$ as $\mathcal W (S, M, \eta) = \tw (\mathcal A)$. The grading structure on $S$ is given by a line field $\eta$, which is a section of the projectivised tangent bundle of $S$. Furthermore, it is shown that, up to Morita equivalence, $\tw (\cala)$ is independent of the full arc system chosen. The category $\cald \calw (S,M,\eta) := \mathrm H^0 \tw (\cala)$ is triangulated. By \cite{HaidenKatzarkovKontsevich} it has as indecomposable objects all (shifts of) curves in $S$ connecting $\gpoint$-points and all  closed curves of winding number zero (with respect to the line field $\eta$) together with an indecomposable local system of finite rank.  It is also shown in \cite{HaidenKatzarkovKontsevich} that $\cald \calw (S,M,\eta)$ is triangle equivalent  to $\ka$,  for any graded gentle algebra $A$ defined by an admissible dissection of $(S,M,\eta)$, see Figure~\ref{fig: HKK disk}. That is, in particular, if $A$ and $B$ are two graded gentle algebras arising from two admissible dissections of $(S,M,\eta)$ then 
 $\ka \simeq \per(B)$. In the language of \cite{HaidenKatzarkovKontsevich} an admissible dissection is called a \textit{full formal arc system}.

 \begin{figure}
    \centering
    \begin{tikzpicture}[x=1em, y=1em]
    \begin{scope}
    \draw[line width=.5pt, overlay] (0,0) circle(3em);
    \draw[fill=yellow!3!green!7,draw=yellow!30!green,line width=.6pt] (225:3em) to[bend right=20] (135:3em) to[bend right=20] (45:3em) to[bend right=20] (-45:3em) to[bend right=20] (-135:3em);
    \node[circle,fill=white,draw=black,inner sep=0,outer sep=0,minimum size=.35em] (BL) at (225:3em) {};
    \node[circle,fill=white,draw=black,inner sep=0,outer sep=0,minimum size=.35em] (BR) at (315:3em) {};
    \node[circle,fill=white,draw=black,inner sep=0,outer sep=0,minimum size=.35em] (TL) at (135:3em) {};
    \node[circle,fill=white,draw=black,inner sep=0,outer sep=0,minimum size=.35em] (TR) at (45:3em) {};
    \foreach \a in {0,90,180,270} {
        \node[circle,fill=green!10!red,inner sep=0,outer sep=0,minimum size=.4em, overlay] at (\a:3em) {};
    }
    \node[font=\scriptsize] at (0,0) {$\mu_4$};
    \node[align=center,anchor=north,font=\footnotesize] at (0,-4em) {$\mathcal A$ from a disk \\ with four marked points \\ and nontrivial $\mu_4$};
    \draw[|-stealth] (4.75em,0) -- (7.25em,0);
    \node[font=\footnotesize] at (6em,.6em) {$\tw$\strut};
    \end{scope}
    \begin{scope}[xshift=12em]
    \draw[line width=.5pt, overlay] (0,0) circle(3em);
    \draw[fill=yellow!3!green!7,draw=yellow!30!green,line width=.6pt] (225:3em) to[bend right=20] (135:3em) to (-45:3em) to[bend right=20] (225:3em);
    \draw[draw=yellow!30!green,line width=.6pt] (45:3em) to (225:3em) (135:3em) to[bend right=20] (45:3em) to[bend right=20] (-45:3em);
    \node[circle,fill=white,draw=black,inner sep=0,outer sep=0,minimum size=.35em] (BL) at (225:3em) {};
    \node[circle,fill=white,draw=black,inner sep=0,outer sep=0,minimum size=.35em] (BR) at (315:3em) {};
    \node[circle,fill=white,draw=black,inner sep=0,outer sep=0,minimum size=.35em] (TL) at (135:3em) {};
    \node[circle,fill=white,draw=black,inner sep=0,outer sep=0,minimum size=.35em] (TR) at (45:3em) {};
    \foreach \a in {0,90,180,270} {
        \node[circle,fill=green!10!red,inner sep=0,outer sep=0,minimum size=.4em, overlay] at (\a:3em) {};
    }
    \node[font=\scriptsize] at (-.4em,-1.3em) {$\mu_3$};
    \node[align=center,anchor=north,font=\footnotesize] at (0,-4em) {$\tw (\mathcal A)$ closure of $\mathcal A$ under \\ shifts and mapping cones \\ with nontrivial $\mu_2, \mu_3, \mu_4$};
    \draw[|-stealth] (4.75em,0) -- (7.25em,0);
    \node[font=\footnotesize] at (6em,.6em) {$\mathrm H^0$\strut};
    \end{scope}
    \begin{scope}[xshift=24em]
    \draw[line width=.5pt, overlay] (0,0) circle(3em);
    \draw[draw=yellow!30!green,line width=.6pt] (225:3em) to[bend right=20] (135:3em) to[bend right=20] (45:3em) to[bend right=20] (-45:3em) to[bend right=20] (-135:3em) to (45:3em) (-45:3em) to (135:3em);
    \node[circle,fill=white,draw=black,inner sep=0,outer sep=0,minimum size=.35em] (BL) at (225:3em) {};
    \node[circle,fill=white,draw=black,inner sep=0,outer sep=0,minimum size=.35em] (BR) at (315:3em) {};
    \node[circle,fill=white,draw=black,inner sep=0,outer sep=0,minimum size=.35em] (TL) at (135:3em) {};
    \node[circle,fill=white,draw=black,inner sep=0,outer sep=0,minimum size=.35em] (TR) at (45:3em) {};
    \foreach \a in {0,90,180,270} {
        \node[circle,fill=green!10!red,inner sep=0,outer sep=0,minimum size=.4em, overlay] at (\a:3em) {};
    }
    \node[align=center,anchor=north,font=\footnotesize] at (0,-4em) {{$\smash{\mathrm H^0} \tw (\mathcal A)$} \\ triangulated category \\ (with trivial $\mu_{\geq 3}$)};
    \end{scope}
    \begin{scope}[xshift=35em]
    \draw[line width=.5pt, overlay] (0,0) circle(3em);
    \draw[draw=yellow!30!green,line width=.6pt] (225:3em) to[bend right=20] (135:3em) (45:3em) to [bend right=20] (-45:3em) to [bend right=20] (-135:3em);
    \node[font=\tiny, color=black!30] at (180:2.1em) {$1$};
    \node[font=\tiny, color=black!30] at (270:2.1em) {$2$};
    \node[font=\tiny, color=black!30] at (0:2.1em) {$3$};
    \draw[-stealth, line width=.5pt, color=black!30] (225:3em) ++ (70:1.2em) arc[start angle=70, end angle=20, radius=1.2em];
    \draw[-stealth, line width=.5pt, color=black!30] (-45:3em) ++ (160:1.2em) arc[start angle=160, end angle=110, radius=1.2em];
    \node[font=\scriptsize, color=black!30] at (225:1.35em) {$a$};
    \node[font=\scriptsize, color=black!30] at (-45:1.35em) {$b$};
    \node[circle,fill=white,draw=black,inner sep=0,outer sep=0,minimum size=.35em] (BL) at (225:3em) {};
    \node[circle,fill=white,draw=black,inner sep=0,outer sep=0,minimum size=.35em] (BR) at (315:3em) {};
    \node[circle,fill=white,draw=black,inner sep=0,outer sep=0,minimum size=.35em] (TL) at (135:3em) {};
    \node[circle,fill=white,draw=black,inner sep=0,outer sep=0,minimum size=.35em] (TR) at (45:3em) {};
    \foreach \a in {0,90,180,270} {
        \node[circle,fill=green!10!red,inner sep=0,outer sep=0,minimum size=.4em, overlay] at (\a:3em) {};
    }
    \node[align=center,anchor=north,font=\footnotesize] at (0,-4em) {$\per (A)$, \,$A \colon 1 \tikztoarg{a} \rel 2 \tikztoarg{b} 3$ \\ perfect derived category \\ from admissible dissection};
    \node at (-5.5em,0) {$\simeq$};
    \end{scope}

    \end{tikzpicture}
    \caption{ Formal and non-formal generators of $\tw (\cala) = \calw(S,M,\eta)$.}
     \label{fig: HKK disk}
\end{figure}

 In \cite{LekiliPolishchuk} an inverse construction is given. Starting with a (homologically smooth) graded gentle algebra $A$, a ribbon graph (which in the terminology of Section \ref{sec: OPS Model} corresponds to the dissection $\zD^*$)  is constructed and from this a marked surface $(S,M)$. Furthermore, from the grading of $A$ a line field $\eta$ on $(S,M)$ is defined such that $\ka$ is triangle equivalent to $\cald \calw (S,M,\eta)$.

Working at the level of pre-triangulated $A_\infty$-categories, Kontsevich conjectured a cosheaf-theoretic description of fully wrapped Fukaya categories of Stein manifolds \cite{Kontsevich}. In \cite{HaidenKatzarkovKontsevich} the authors show that such a cosheaf-theoretical construction holds for partially wrapped Fukaya categories of surfaces with stops and relate it to graded gentle algebras. In the terminology of Section~\ref{sec: OPS Model}, they show that given a graded marked surface $(S,M,\eta)$ with dual admissible dissections $\zD$ and $\zD^*$, the category $\tw (\cala)$ is the category of global sections of a cosheaf of pre-triangulated $A_\infty$-categories on $\zD^*$ (which plays the role of a Lagrangian core of the surface).

This cosheaf-theoretical point of view also sheds light on gentle algebras themselves. Namely, a graded gentle algebra $A$ can be viewed as a homotopy colimit of building blocks consisting of graded quivers of type $\mathbb A$ and $\widetilde{\mathbb A}$ modulo the ideal of relations generated by all paths of length two, where the colimit is taken in the homotopy category of dg categories in the sense of \cite{KellerDG,Tabuada,Toen}, see Figure~\ref{fig: Cosheaf Orbifold} and its caption.

We further note that in \cite{HaidenKatzarkovKontsevich}, the surface $(S, M, \eta)$ arises as the metric completion of a flat metric on a Riemann surface $\Sigma$ with exponential-type singularities. It gives rise to an isomorphism between the space of Bridgeland stability conditions on $\mathcal D \mathcal W (S, M, \eta)$ and the space of flat structures $\mathcal M (\Sigma)$ on $\Sigma$ \cite{HaidenKatzarkovKontsevich,Takeda}.

\section{Exceptional sequences: Existence and  braid group action}

Exceptional sequences are a tool to decompose a triangulated category into smaller pieces and a full exceptional sequence can be viewed as a particular generator of a triangulated category. They first arose in the work of Beilinson in the study of $\D^b(\operatorname{coh} \mathbb{P}^n)$ where he showed that the sequence $(\calo, \calo(1), \ldots,  \calo(n))$ of line bundles forms a full exceptional sequence. Generalisations to other triangulated categories quickly developed, with full exceptional sequences emerging as one of the bridges between algebraic geometry and representation theory of finite dimensional algebras \cite{Bondal}. 

A sequence $\mathbf{E}=(E_1, \ldots, E_n)$ in a triangulated category $\calt$ is \textit{exceptional}, if $\Hom_\calt(E_i, E_i[\mathbb Z_{\neq0}])= 0$,  $\Hom_\calt(E_i, E_i)$ is a division ring and $\Hom_\calt(E_i, E_j[\mathbb{Z}]) = 0$, for $j >i$. An exceptional sequence $\mathbf{E}$ is \textit{full} if $\thick_\calt(\mathbf{E})=\calt$. 
We note that shifting the $E_i$ in an exceptional sequence gives rise to another exceptional sequence, so there is a natural action of $\mathbb{Z}^n$ on the set of full exceptional sequences in $\calt$.

\subsection{Exceptional sequences in the derived category of a graded gentle algebra}

In general, it is very difficult to find full exceptional sequences in a triangulated category or indeed to determine whether there exist any at all. In the case of the perfect derived category $\ka$, of a graded gentle algebra $A$ with associated graded marked surface $(S,M,\eta)$, it is shown in   \cite{ChangJinSchroll, ChangSchrollexcp}  that full exceptional sequences exist if and only if there are no $\rpoint$-points in the interior of $S$ (equivalently, $A$ is homologically smooth) and if there are at least two $\gpoint$-points in the boundary. Moreover, if they exist, full exceptional sequences correspond bijectively to those admissible surface dissections where the associated quiver has no oriented cycle \cite{ChangJinSchroll,ChangSchrollexcp}. This is a criterion which can be formulated solely in terms of the surface, see \cite{ChangSchrollexcp}. The characterisation of full exceptional sequences in terms of admissible surface dissections shows that in case $\per (A)$ admits a full exceptional sequence, any exceptional sequence of the same length is full, confirming \cite[Conjecture 1.10]{Kuznetsov} for $\per (A)$.
Exceptional sequences are closely related to semiorthogonal decompositions which for the derived categories of graded gentle algebras have been studied in terms of the associated surface in \cite{KoprivaStovicek}.

\subsection{Example of non-transitivity of the braid group action on full exceptional sequences}

By \cite{GorodentsevRudakov} there is a natural mutation operation on exceptional sequences which induces an action of the braid group $\frak{B}_n$ on the set of full exceptional sequences in $\calt$ \cite{Bondal}. In \cite{BondalPolishchuk} this action was further studied and  it was conjectured that the action of $\Z^n \rtimes \mathfrak{B}_n$ should be transitive. This has been shown to hold in many cases, for example, for derived categories of hereditary algebras \cite{CBexceptional, RingelExceptional}, for derived categories of coherent sheaves over del Pezzo surfaces, projective planes and weighted projective lines \cite{KussinMeltzer, Meltzer}, and for the perfect derived category of a gentle algebra where the associated surface is of genus zero \cite{ChangSchrollexcp}. 

Using Birman--Hilden theory and Hurwitz systems in the theory of branched coverings, it is shown in \cite{ChangHaidenSchroll} that for certain graded gentle algebras the action of $\mathbb{Z}^n \rtimes \frak{B}_n$ on the set of full exceptional sequences of $\ka$ is not only not transitive, but decomposes into infinitely many orbits.  This applies, in particular, to the class of graded gentle algebras whose associated graded marked surface $(S, M, \eta)$ has exactly two marked points and genus greater than 1.   These are the graded gentle algebras derived equivalent to graded gentle algebras whose  quivers, for $n \geq 3$,  are given by: 
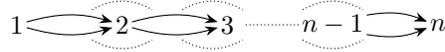
\[
\begin{tikzpicture}[>=Stealth, baseline=(current bounding box.center)]
\tikzset{every node/.style={font=\small, inner sep=1pt}}

% 3) 1 ⇉ 2 ⇉ 3
\begin{scope}[shift={(0.0,0)}]
  \node (c1) at (0,0) {$1$};
  \node (c2) at (1.4,0) {$2$};
  \node (c3) at (2.8,0) {$3$};
  \node (d1) at  (3.0,0) {}; 
  \node (d2) at (3.8,0) {}; 
  \node (c4) at (4.2,0) {$n-1$};
  \node (c5) at (5.6,0) {$n$};
  
  % double arrows 1 -> 2
  \draw[->,overlay] (c1) to[bend left=16,overlay]  (c2);
  \draw[->,overlay] (c1) to[bend right=16,overlay] (c2);
  \draw[densely dotted,overlay] (1,0.2) to [bend left,overlay] (1.8,0.2);
   \draw[densely dotted,overlay] (1,-0.2) to [bend right,overlay] (1.8,-0.2);
  % double arrows 2 -> 3
  \draw[->] (c2) to[bend left=16,overlay]  (c3);
  \draw[->,overlay] (c2) to[bend right=16] (c3);
 \draw[densely dotted,overlay] (2.2,0.2) to [bend left,overlay] (3.0,0.2);
   \draw[densely dotted,overlay] (2.2,-0.2) to [bend right,overlay] (3.0,-0.2);

   \draw[densely dotted] (d1) to (d2);

  % double arrows n-1 -> n
    \draw[->,overlay] (c4) to[bend left=16,overlay]  (c5);
    \draw[->,overlay] (c4) to[bend right=16,overlay] (c5);
    \draw[densely dotted,overlay] (3.8,0.2) to [bend left,overlay] (4.6,0.2);
   \draw[densely dotted] (3.8,-0.2) to [bend right] (4.6,-0.2);
\end{scope}

\end{tikzpicture}
\]

A new algebraic proof of the non-transitivity of the braid group action for these algebras was recently given in  \cite{NakagoTakahashi}.

\section{Derived equivalences}

In this section we show how  the surface model of Section~\ref{sec: OPS Model} can be used to obtain results on derived equivalences and give a complete derived invariant for graded gentle algebras. 

%{\color{magenta} 
%Using the surface model of Section \ref{sec: OPS Model}, it is also possible to produce derived equivalences and to give a complete combinatorial derived invariant in terms of winding numbers. }

\subsection{Derived equivalences: Silting and tilting objects as admissible surface dissections} Another type of generator in a triangulated category  is given by  the so-called \textit{silting objects} introduced in \cite{KellerVossieck}.   By \cite{AiharaIyama} every full exceptional sequence gives rise to a silting object. This has the immediate consequence that whenever $\calt$ admits  a full exceptional sequence, it also has a silting object. Just as full exceptional sequences do not always exist, silting objects need not exist either.

Let $\calt$ be a triangulated category. An object $X \in \calt$ is called \textit{pre-silting} if  $\Hom_\calt(X, X[\mathbb{Z}_{>0}])= 0$.  If in addition $\thick_\calt(X)= \calt$, then $X$ is called \textit{silting}. 
We note that, for a non-positively graded  dg algebra $A$, the algebra $A$ itself is a silting object in $\ka$.
Silting objects give rise to derived equivalences, see, for example, \cite{KellerDerivingDG,KellerConstrofTriaEq}. 

Silting objects are a generalisation of tilting objects, where in the context of associative algebras, the notion of  a  tilting complex in the derived category as a generalisation of (generalised) tilting modules has been introduced in \cite{Rickard89} in the study of derived Morita equivalences, see also \cite{KellerConstrofTriaEq}. For an expository  account of tilting theory, see, for example,   \cite{HandBookOfTiltingTheory} or \cite{JassoKrauseSchroll}. Given a triangulated category $\calt$, an object $T \in \calt$ is called a \textit{tilting object}
if $\Hom_\calt(T, T[\mathbb Z_{\neq0}])=0$ and if $\thick_\calt(T) \simeq \calt$. For example, for a finite dimensional algebra $A$, if $\mathcal T = \ka$ and $T$ a tilting object in $\calt$ with  $B = \End_{\calt}(T)$ then $\ka $ and $\per(B)$ as well as $\dba$ and $\mathcal{D}^b(B)$ are triangle equivalent. 
In this setting, any derived equivalence is  induced by a tilting complex.  

 A key technique for constructing derived equivalences is to construct tilting and silting objects and to calculate their endomorphism algebras. For  graded gentle algebras,  the associated graded marked surface is a helpful tool for this. More specifically, in 
\cite{AmiotPlamondonSchroll, Opper}  it is shown  that for zero-graded gentle algebras tilting complexes correspond to  certain admissible surface dissections. In \cite{JinSchrollWang} it is shown that the same holds  for graded gentle algebras, in that  silting objects also correspond to admissible surface dissections in this case. Note, however, that there is no geometric criterion to determine  whether an  admissible dissection  corresponds to a silting or tilting object. 

\subsection{Existence of silting objects} While in the case of a zero-graded algebra, tilting (and hence silting) objects always exist (the algebra itself is a tilting object), in the case of a graded algebra that is not  non-positively graded, silting objects might not exist in $\ka$. 
Indeed, using reduction techniques given by the correspondence of surface cuts and recollements of derived categories in \cite{ChangJinSchroll},  in \cite{JinSchrollWang} we give a complete characterisation of when a  homologically smooth (not necessarily proper) graded gentle algebra   has silting objects. That is, such an algebra $A$ does not admit a silting object if and only if its quiver has an oriented cycle $C=a_m \cdots a_1$  without relations and $|C| = |a_1| + \dotsb + |a_m| > 0$, or if $A$ is isomorphic to the graded gentle algebra $A'$ with quiver 
%\begin{center}
%\begin{figure}[ht]
 %   \centering
 %  \begin{tikzpicture}[>=stealth, width=3cm,]
 \begin{center}
 \begin{tikzpicture}[>=stealth, node distance=2cm]  
  % Nodes
  \node (1) {1};
  \node[right of=1] (2) {2};

  % Arrows with labels
  \draw[->,font=\scriptsize,overlay] (1) to[bend left=30] node[overlay,above=-.3ex] {$a$} (2);   % top arrow 1 -> 2
  \draw[->,font=\scriptsize] (2) to node[midway,above=-.3ex] {$b$} (1);          % middle arrow 2 -> 1
  \draw[->,font=\scriptsize,overlay] (1) to[bend right=30] node[above=-.3ex,overlay] {$c$} (2);  % bottom arrow 1 -> 2
\end{tikzpicture}
\end{center}
and relations $ba, cb$ and grading $|a|+|b| = 1 = |b| + |c| $. 
Geometrically, $A$ has a silting object if and only if the associated graded surface $S$ has no $\gpoint$-points in the interior with negative winding number, or if $S$ is the torus with one boundary component with exactly one $\gpoint$-point on that boundary, such that the winding number of each non-separating simple closed curve is zero. The observation that the algebra $A'$ has no silting object already appeared in \cite{ChangJinSchroll} and was also proved in \cite{LiuZhou}. 

Note that, for the algebra $A'$  (which is such that $\per(A')$ admits no silting object), it has recently been shown in \cite{HaidenWu} that   $\per (A')$ does also not admit any Bridgeland stability conditions.

\subsection{Not every pre-silting object can be completed to a silting object}\label{sec:presilting vs partial silting}

The first example of a finite dimensional algebra that has a pre-silting object that cannot be completed to a silting object appeared in \cite{LiuZhou}. In \cite{JinSchrollWang}, building on \cite{ChangJinSchroll}, it is then shown that infinitely many finite dimensional gentle algebras have pre-silting objects which cannot be completed to a silting object. More precisely,  it is shown that, a finite dimensional (zero-graded) gentle algebra $A$ of finite global dimension has a pre-silting object that cannot be completed to a silting object, as soon as the genus of the associated graded surface is strictly greater than 1.

\subsection{Complete derived invariant for graded gentle algebras}

The surface model gives insight into the derived equivalences of graded gentle algebras. It follows from the construction in \cite{HaidenKatzarkovKontsevich} that all graded gentle algebras that arise from surface dissections of a graded marked surface are derived equivalent, see also Section~\ref{sec: Fukaya}.

On the other hand, it is not clear, a priori, whether graded gentle algebras which are derived equivalent have homeomorphic surface models. This turns out to be the case, as has been shown for zero-graded gentle algebras independently in   \cite{AmiotPlamondonSchroll} and \cite{Opper} using tilting theory, and for (homologically smooth not necessarily proper) graded gentle algebras in \cite{JinSchrollWang} using silting theory: Given two graded gentle algebras $A$ and $A'$ with surface models $(S, M, \eta)$ and $(S', M', \eta')$, respectively, the algebras $A$ and $A'$ are derived equivalent if and only if there exists an orientation preserving homeomorphism $\phi \colon S \to S'$ such that $\phi (M) = M'$ and $\phi_* (\eta)$ is homotopic to $\eta'$.

In other words, the graded marked surface can be viewed as a complete derived invariant for graded gentle algebras \cite{AmiotPlamondonSchroll,JinSchrollWang, Opper}. A purely combinatorial derived invariant $\mathrm{AG} (A)$ of a zero-graded gentle algebra $A$ was already introduced in \cite{AvellaGeiss}. It is independently shown in \cite{LekiliPolishchuk} and \cite{OpperPlamondonSchroll} that $\mathrm{AG} (A)$ has a geometric interpretation. It encodes the number of boundary components, the number of marked points on each boundary component and the winding numbers of the boundary components. In particular, this result extends the invariant $\mathrm{AG}$ to graded gentle algebras.
However, it was known from \cite{Amiot} that this is not a complete invariant. In other words, $\mathrm{AG} (A) = \mathrm{AG} (A')$ does not necessarily imply that $A$ and $A'$ are derived equivalent. In \cite{Amiot, AmiotGrimeland} this invariant was extended to a complete derived invariant  for the subclass of gentle algebras defined in \cite{DavidSchiffler}. More generally, for arbitrary graded gentle algebras, the invariant $\mathrm{AG} (A)$ was augmented in  \cite{LekiliPolishchuk} with two further invariants involving the winding numbers of a set of generators of the fundamental group of the associated closed surface. This augmented invariant was shown to be a complete derived invariant in \cite{AmiotPlamondonSchroll, Opper}
for the zero-graded case and in \cite{JinSchrollWang} for the graded case. Very recently another proof for the graded case has appeared in \cite{Opper26}.

\section{$A_\infty$-deformations of  gentle algebras and topological Fukaya categories of orbifold surfaces}

In this section, we discuss how the Hochschild cohomology (and  the Tamarkin--Tsygan-calculus) of a graded gentle algebra is encoded on the associated graded marked surface. In particular, the knowledge of the second Hochschild cohomology groups for graded gentle algebras enables the  study of their $A_\infty$-deformations.
The construction of topological Fukaya categories of graded orbifold surfaces then makes it possible to show that $A_\infty$-deformations of graded gentle algebras correspond to partial compactifications of the associated graded marked surfaces.

\subsection{Hochschild cohomology and the Tamarkin--Tsygan calculus for graded gentle algebras}\label{sec: HH}

The notion of a Tamarkin--Tsygan calculus (TT-calculus) or differential calculus was introduced in \cite{TamarkinTsygan}. The first class of examples of such a calculus is given by the Hochschild cohomology and homology of associative algebras with their associated structures induced by the cup and cap products, the Gerstenhaber bracket and the Connes differential  \cite{GelfandDaletskiiTsygan, TamarkinTsygan}. 
The TT-calculus of an algebra is invariant under derived equivalence, as has been shown in a series of papers by different authors concluding with \cite{ArmentaKellerTT}, see also the references within.

\begin{figure}
    \centering
    \begin{tikzpicture}[scale=1.1, x=1em, y=1em]
    \begin{scope}
        \draw[line width=.3pt] (-4em,-2.2em) to ++ (0, .1em);
        \draw[line width=.3pt] (33.05em,-2.2em) to ++ (0, .1em);
        \draw[line width=.4pt, overlay] (-4em, 2.8em) rectangle ++(37.05em, -5.6em);
        \draw[line width=.4pt, overlay] (10.2em, 2.8em) to[overlay] ++(0, -5.6em) (21.8em, 2.8) to[overlay] ++(0, -5.6em);
        \draw[line width=.5pt] (0,0) circle(1em);
        \draw[draw=black!20, line width=.3pt] (0,0) circle(1.5em);
        \node[font=\scriptsize, color=black!20] at (-3:1.9em) {$C$};
        \node[circle,fill=green!10!red,inner sep=0,outer sep=0,minimum size=.4em] at (90:1em) {};
        \draw[draw=green!10!red, line width=.6pt] (190:3em) to[out=45, in=150] (90:1em) (135:2.5em) to[out=0, in=120] (90:1em) (-10:3em) to[in=30, out=135] (90:1em) (45:2.5em) to[out=180, in=60] (90:1em);
        \draw[draw=yellow!30!green, line width=.6pt] (-190:3em) to[out=-45, in=-150] (-90:1em) (-135:2.5em) to[out=0, in=-120] (-90:1em) (10:3em) to[in=-30, out=-135] (-90:1em) (-45:2.5em) to[out=-180, in=-60] (-90:1em);
        \begin{scope}[yshift=1em]
        \draw[-stealth, overlay] (147:1.25em) arc[start angle=147, end angle=175, radius=1.25em];
        \draw[-stealth, overlay] (4:1.25em) arc[start angle=4, end angle=32, radius=1.25em];
        \draw[-stealth, overlay] (36:1.25em) arc[start angle=36, end angle=64, radius=1.25em];
        \node[font=\scriptsize, overlay] at (78:1.25em) {.};
        \node[font=\scriptsize, overlay] at (90:1.25em) {.};
        \node[font=\scriptsize, overlay] at (102:1.25em) {.};
        \draw[-stealth, overlay] (116:1.25em) arc[start angle=116, end angle=144, radius=1.25em];
        \node[font=\scriptsize, right] at (1.15em,.3em) {$a_1$};
        \node[font=\scriptsize, left] at (-1.15em,.4em) {$a_m$};
        \end{scope}
        \begin{scope}[yshift=-1em]
        \draw[-stealth, overlay] (-147:1.25em) arc[start angle=-147, end angle=-175, radius=1.25em];
        \draw[-stealth, overlay] (-4:1.25em) arc[start angle=-4, end angle=-32, radius=1.25em];
        \draw[-stealth, overlay] (-36:1.25em) arc[start angle=-36, end angle=-64, radius=1.25em];
        \node[font=\scriptsize, overlay] at (-78:1.25em) {.};
        \node[font=\scriptsize, overlay] at (-90:1.25em) {.};
        \node[font=\scriptsize, overlay] at (-102:1.25em) {.};
        \draw[-stealth, overlay] (-116:1.25em) arc[start angle=-116, end angle=-144, radius=1.25em];
        \node[font=\scriptsize, right] at (1.15em,-.3em) {$b_1$};
        \node[font=\scriptsize, left] at (-1.1em,-.25em) {$b_n$};
        \end{scope}
        \node[circle,fill=white,draw=black,inner sep=0,outer sep=0,minimum size=.35em] at (-90:1em) {};
        \begin{scope}[xshift=6em]
        \node[font=\footnotesize] at (0,1em) {$x_{a,b}$}; 
        \node[font=\footnotesize] at (0,0) {\rotatebox{-90}{$\in$}};
        \node[font=\footnotesize] at (0,-1.2em) {$\HH^{\mathrm w_\eta (C) + 1} (A)$};
        \end{scope}
        \begin{scope}[xshift=13em] % middle
        \draw[draw=black!20, line width=.3pt] (0,0) circle(1.35em);
        \node[font=\scriptsize, color=black!20] at (225:1.85em) {$C$};
        \node[circle,fill=green!10!red,inner sep=0,outer sep=0,minimum size=.4em] at (0,0) {};
        \draw[draw=green!10!red, line width=.6pt] (0,0) to (20:1.5em);
        \draw[draw=green!10!red, line width=.6pt] (0,0) to (85:1.5em);
        \draw[draw=green!10!red, line width=.6pt] (0,0) to (150:1.5em);
        \draw[draw=green!10!red, line width=.6pt] (0,0) to (-45:1.5em);
        \foreach \a in {-110, -45, 20, 85, 150} {
        \draw[-stealth] (\a+2:1em) arc[start angle=\a+2, end angle=\a+61, radius=1em];
        }
        \node[font=\scriptsize] at (51.5+180:1em) {.};
        \node[font=\scriptsize] at (51.5+180+12:1em) {.};
        \node[font=\scriptsize] at (51.5+180-12:1em) {.};
        \node[font=\scriptsize] at (20-32.5:1.65em) {$a_m$};
        \node[font=\scriptsize] at (85-32.5-5:1.55em) {$a_1$};
        \node[font=\scriptsize] at (150-32.5:1.5em) {$a_2$};
        \begin{scope}[xshift=5em]
        \node[font=\footnotesize] at (0,1em) {$x_{a^i\!,\,\mathrm s (a)}$}; 
        \node[font=\footnotesize] at (0,0) {\rotatebox{-90}{$\in$}};
        \node[font=\footnotesize] at (0,-1.2em) {$\HH^{i\, \mathrm w_\eta (C)} (A)$};
        \end{scope}
        \begin{scope}[xshift=11.5em]
        \end{scope}
        \end{scope}
        \begin{scope}[xshift=24.8em] % right
        \draw[draw=black!20, line width=.3pt] (0,0) circle(1.35em);
        \node[font=\scriptsize, color=black!20] at (-45:1.8em) {$C$};
        \draw[draw=yellow!30!green, line width=.6pt] (0,0) to (180-20:1.5em);
        \draw[draw=yellow!30!green, line width=.6pt] (0,0) to (180-85:1.5em);
        \draw[draw=yellow!30!green, line width=.6pt] (0,0) to (180-150:1.5em);
        \draw[draw=yellow!30!green, line width=.6pt] (0,0) to (180+45:1.5em);
        \foreach \a in {110, 45, -20, -85, -150} {
        \draw[-stealth] (180+\a-2:1em) arc[start angle=180+\a-2, end angle=180+\a-61, radius=1em];
        }
        \node[font=\scriptsize] at (-51.5:1em) {.};
        \node[font=\scriptsize] at (-51.5+12:1em) {.};
        \node[font=\scriptsize] at (-51.5-12:1em) {.};
        \node[font=\scriptsize] at (180-20+32.5:1.5em) {$b_n$};
        \node[font=\scriptsize] at (180-85+32.5-5:1.55em) {$b_1$};
        \node[font=\scriptsize] at (180-150+32.5:1.5em) {$b_2$};
        \node[circle,fill=white,draw=black,inner sep=0,outer sep=0,minimum size=.35em] at (0,0) {};
        \begin{scope}[xshift=5em]
        \node[font=\footnotesize] at (0,1em) {$x_{\mathrm s(b){\,\!,\, b^i}}$}; 
        \node[font=\footnotesize] at (0,0) {\rotatebox{-90}{$\in$}};
        \node[font=\footnotesize] at (0,-1.2em) {$\HH^{i\, \mathrm w_\eta (C)} (A)$};
        \end{scope}
        \end{scope}
    \end{scope}
    \end{tikzpicture}
    \caption{Geometric representation of the three types of non-Eulerian Hochschild cocycles generating $\HH^*(A)$ as a graded commutative algebra. For $ a = a_m \dotsb a_1$ and $b = b_n \dotsb b_1$, the generator $x_{a,b}$ corresponds to a morphism $A^{\otimes m} \to A$ which maps $a_m \otimes \cdots \otimes a_1$ to $b$ (left), $x_{a^i, \mathrm s (a)} \colon A^{\otimes im} \to A$ maps $(a_{j+m-1} \otimes \dotsb  \otimes  a_j)^{ \otimes i}$ to $e_{\mathrm s (a_j)}$ for all $0 \leq j < m$ (middle) and, dually, $x_{\mathrm s (b), b^i} \colon \Bbbk Q_0 \to A$ maps $e_{\mathrm s (b_j)}$ to $(b_{j+n-1} \dotsb  b_j)^i$ for all $0 \leq j < n$ (right), indices being labelled cyclically, and $i = 2$ unless $\mathrm w_\eta (C)$ is even or $\mathrm{char} (\Bbbk) = 2$ in which case $i = 1$.}  
     \label{fig: Hochschild}
\end{figure}
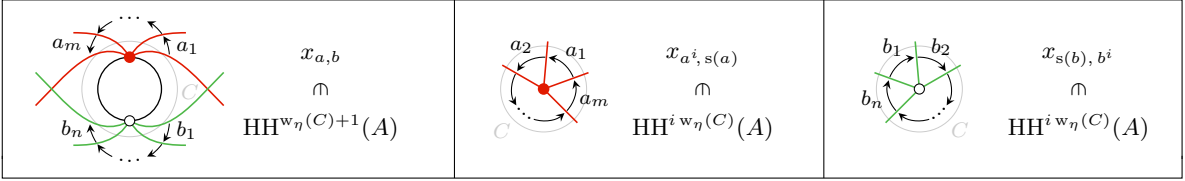

As the calculation of all the ingredients involved in the TT-calculus is generally very difficult, there have been very few calculations of the complete calculus. One of the rare cases, where a complete calculation is possible is  the class of gentle algebras. The complete TT-calculus of a (possibly infinite dimensional) zero-graded gentle algebra $A = \Bbbk Q/I$ is calculated in \cite{ChaparroSchrollSolotarSuarez}. In particular, it is shown that the generators of the Hochschild cohomology as a graded commutative algebra 
\cite{ChaparroSchrollSolotarSuarez}
have a geometric interpretation in terms of the surface model of the bounded derived category. The so-called Eulerian derivations in $\HH^1(A)$ are in bijection with a set of generators of the fundamental group of the surface with boundary (which by \cite{OpperPlamondonSchroll} is isomorphic to the fundamental group of the quiver of $A$).
All other generators of $\HH^* (A)$ are in bijection with the set consisting of the boundary components with one marked point as well as with $\rpoint$-points and $\gpoint$-points in the interior of the surface $(S,M,\eta)$ associated to $A$. Furthermore, the winding number around the boundary or $\rpoint$- or $\gpoint$-points is related to the degree in which the corresponding generator of $\HH^*(A)$ lies, see Figure~\ref{fig: Hochschild}. 

It follows from the complete description of the structure of $ (\HH^*(A),\smile)$ as a graded commutative algebra in terms of generators and relations in 
\cite{ChaparroSchrollSolotarSuarez}, that it is a quadratic algebra. In many cases  it even is 
quadratic monomial and hence Koszul.
This holds, for example, if $A$ is both finite dimensional and of finite global dimension but it also holds for many infinite dimensional gentle algebras of infinite global dimension, see 
\cite{ChaparroSchrollSolotarSuarez}.

Using a graded version of the Bardzell resolution \cite{Bardzell}, these calculations have been extended in \cite{BarmeierSchrollWangDef} to calculate $\HH^2(A)$ for a graded gentle algebra $A$. This is then used to study the $A_\infty$-deformations of graded gentle algebras. 
In particular, it is shown that a basis of $\HH^2(A)$ is in bijection with the union of 
\newline $\bullet$ the set of interior $\rpoint$-points of winding number 1 or 2,
\newline $\bullet$ the set of boundary components with exactly one $\gpoint$-point and winding number 1
\newline $\bullet$ the set of interior $\gpoint$-points of winding number 1 or 2.

In  \cite{BianSchrollSolotarWangWen}, it is shown that the whole of the graded commutative structure and indeed the whole of the Gerstenhaber structure of $\HH^*(A)$ of a graded (skew-)gentle algebra can be described in terms of the associated graded marked (orbifold) surface in exactly the same way as illustrated in Figure~\ref{fig: Hochschild}. A calculation of the Gerstenhaber structure of the Hochschild cohomology in the graded gentle case was also carried out in \cite{Opper26}. 

\subsection{$A_\infty$-deformations of graded gentle algebras}

The knowledge of the Hochschild cohomology can be used to describe the $A_\infty$-deformation theory of graded gentle algebras. In \cite{BarmeierSchrollWangDef}, using the invariance of the deformation theory of dg categories under derived equivalence \cite{KellerHH}, see also \cite{LowenVandenBergh}, an explicit family $\{ \per (A_\lambda) \}_{\lambda \in \Bbbk^d}$ of deformations of the  (dg enhanced) perfect derived category of a graded gentle algebra  $A$ is constructed. Each $A_\lambda$ can be viewed as a dg deformation of the graded gentle algebra $A = A_0$. As the Hochschild $2$-cocycles of $A$  naturally correspond to boundary components with winding number $1$ or $2$ (see Section \ref{sec: HH}) a natural question is, whether the corresponding deformations have a natural geometric interpretation as well. This is indeed the case: the individual fibers $\per (A_\lambda)$ of the family of $A_\infty$-deformations correspond to partial compactifications of the graded marked surface of  $A$. Specifically,  the boundary components of winding number $1$ are replaced by disks with orbifold points of order $2$ and the boundary components of winding number $2$ by smooth disks. This reflects the fact that these are precisely the two cases when the line field can be extended smoothly to the disks or orbifold disks that are glued into the boundary components.

These results  give a first account of the relationship between deformations of Fukaya categories and partial compactifications in the context of partially wrapped Fukaya categories, proposed in the fully wrapped case in \cite{SeidelICM} and studied for example in \cite{PerutzSheridan, Sheridan}.

The $A_\infty$-deformations of infinite dimensional $\mathbb Z / 2$-graded gentle algebras are studied in \cite{BocklandtvandeKreeke}
building on \cite{Bocklandt}.

\subsection{Topological Fukaya categories of orbifold surfaces with stops}\label{sec: Fukaya orbifold}

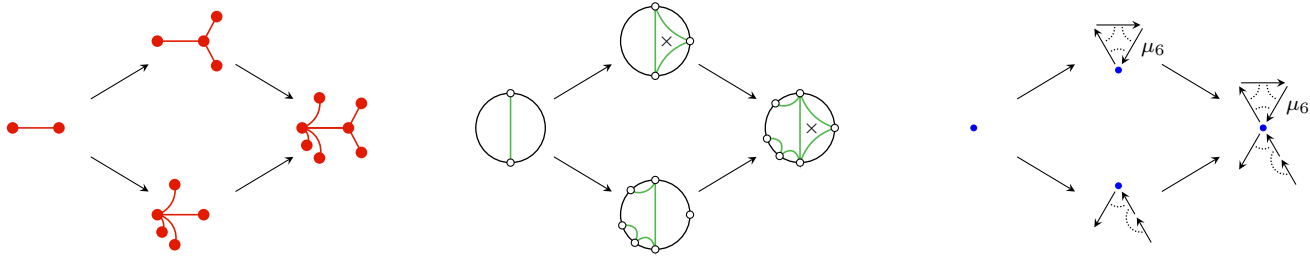
\begin{figure}
    \centering
    \begin{tikzpicture}
        \begin{scope}[xshift=-16em] % left
        \draw[-stealth] (-3.5em,1em) to (-1.5em,2.2em);
        \draw[-stealth] (-3.5em,-1em) to (-1.5em,-2.2em);
        \draw[-stealth] (1.5em,2.2em) to (3.5em,1em);
        \draw[-stealth] (1.5em,-2.2em) to (3.5em,-1em);
        \begin{scope}[xshift=-5em] % L
            \node[circle,fill=green!10!red,inner sep=0,outer sep=0,minimum size=.4em] at (-1.2em,0) {};
            \node[circle,fill=green!10!red,inner sep=0,outer sep=0,minimum size=.4em] at (.4em,0) {};
            \draw[draw=green!10!red, line width=.6pt] (-1.2em,0) to (.4em,0);
        \end{scope}
        \begin{scope}[yshift=3em] % T
            \node[circle,fill=green!10!red,inner sep=0,outer sep=0,minimum size=.4em] at (-1.2em,0) {};
            \node[circle,fill=green!10!red,inner sep=0,outer sep=0,minimum size=.4em] at (.4em,0) {};
            \node[circle,fill=green!10!red,inner sep=0,outer sep=0,minimum size=.4em] at (45:1.2em) {};
            \node[circle,fill=green!10!red,inner sep=0,outer sep=0,minimum size=.4em] at (-45:1.2em) {};
            \draw[draw=green!10!red, line width=.6pt] (-1.2em,0) to (.4em,0) to (45:1.2em) (.4em,0) to (-45:1.2em);
        \end{scope}
        \begin{scope}[yshift=-3em] % B
            \node[circle,fill=green!10!red,inner sep=0,outer sep=0,minimum size=.4em] at (-1.2em,0) {};
            \node[circle,fill=green!10!red,inner sep=0,outer sep=0,minimum size=.4em] at (.4em,0) {};
            \node[circle,fill=green!10!red,inner sep=0,outer sep=0,minimum size=.4em] at (240:1.2em) {};
            \node[circle,fill=green!10!red,inner sep=0,outer sep=0,minimum size=.4em] at (120:1.2em) {};
            \node[circle,fill=green!10!red,inner sep=0,outer sep=0,minimum size=.4em] at (210:1.2em) {};
            \draw[draw=green!10!red, line width=.6pt] (-1.2em,0) to (.4em,0) (-1.2em,0) to[bend left=45] (210:1.2em) (-1.2em,0) to[bend left=50] (240:1.2em) (-1.2em,0) to[bend right=40] (120:1.2em);
        \end{scope}
        \begin{scope}[xshift=5em] % R
            \node[circle,fill=green!10!red,inner sep=0,outer sep=0,minimum size=.4em] at (-1.2em,0) {};
            \node[circle,fill=green!10!red,inner sep=0,outer sep=0,minimum size=.4em] at (.4em,0) {};
            \node[circle,fill=green!10!red,inner sep=0,outer sep=0,minimum size=.4em] at (240:1.2em) {};
            \node[circle,fill=green!10!red,inner sep=0,outer sep=0,minimum size=.4em] at (120:1.2em) {};
            \node[circle,fill=green!10!red,inner sep=0,outer sep=0,minimum size=.4em] at (210:1.2em) {};
            \node[circle,fill=green!10!red,inner sep=0,outer sep=0,minimum size=.4em] at (45:1.2em) {};
            \node[circle,fill=green!10!red,inner sep=0,outer sep=0,minimum size=.4em] at (-45:1.2em) {};
            \draw[draw=green!10!red, line width=.6pt] (-1.2em,0) to (.4em,0);
            \draw[draw=green!10!red, line width=.6pt] (-1.2em,0) to (.4em,0) (-1.2em,0) to[bend left=45] (210:1.2em) (-1.2em,0) to[bend left=50] (240:1.2em) (-1.2em,0) to[bend right=40] (120:1.2em) (.4em,0) to (45:1.2em) (.4em,0) to (-45:1.2em);
        \end{scope}
        \end{scope} % end left
        \begin{scope} % middle
        \draw[-stealth] (-3.5em,1em) to (-1.5em,2.2em);
        \draw[-stealth] (-3.5em,-1em) to (-1.5em,-2.2em);
        \draw[-stealth] (1.5em,2.2em) to (3.5em,1em);
        \draw[-stealth] (1.5em,-2.2em) to (3.5em,-1em);
        \begin{scope}[xshift=-5em] % L
            \draw[line width=.5pt] (0,0) circle(1.2em);
            \draw[line width=.6pt, draw=yellow!30!green] (90:1.2em) to (270:1.2em);
            \node[circle,fill=white,draw=black,inner sep=0,outer sep=0,minimum size=.25em] at (90:1.2em) {};
            \node[circle,fill=white,draw=black,inner sep=0,outer sep=0,minimum size=.25em] at (270:1.2em) {};
        \end{scope}
        \begin{scope}[yshift=3em] % T
            \draw[line width=.5pt] (0,0) circle(1.2em);
            \draw[line width=.6pt, draw=yellow!30!green] (90:1.2em) to (270:1.2em) to[bend left=25] (0:1.2em) to [bend left=25] (90:1.2em);
            \node[circle,fill=white,draw=black,inner sep=0,outer sep=0,minimum size=.25em] at (0:1.2em) {};
            \node[circle,fill=white,draw=black,inner sep=0,outer sep=0,minimum size=.25em] at (90:1.2em) {};
            \node[circle,fill=white,draw=black,inner sep=0,outer sep=0,minimum size=.25em] at (270:1.2em) {};
            \node[font=\tiny] at (0:.4em) {$\times$};
        \end{scope}
        \begin{scope}[yshift=-3em] % B
            \draw[line width=.5pt] (0,0) circle(1.2em);
            \draw[line width=.6pt, draw=yellow!30!green] (90:1.2em) to (270:1.2em) to[bend right=60, looseness=1.6] (234:1.2em) to[bend right=60, looseness=1.6] (198:1.2em) (90:1.2em) to[bend left=60] (135:1.2em);
            \node[circle,fill=white,draw=black,inner sep=0,outer sep=0,minimum size=.25em] at (0:1.2em) {};\node[circle,fill=white,draw=black,inner sep=0,outer sep=0,minimum size=.25em] at (90:1.2em) {};
            \node[circle,fill=white,draw=black,inner sep=0,outer sep=0,minimum size=.25em] at (135:1.2em) {};
            \node[circle,fill=white,draw=black,inner sep=0,outer sep=0,minimum size=.25em] at (198:1.2em) {};
            \node[circle,fill=white,draw=black,inner sep=0,outer sep=0,minimum size=.25em] at (234:1.2em) {};
            \node[circle,fill=white,draw=black,inner sep=0,outer sep=0,minimum size=.25em] at (270:1.2em) {};
        \end{scope}
        \begin{scope}[xshift=5em] % R
            \draw[line width=.5pt] (0,0) circle(1.2em);
            \draw[line width=.6pt, draw=yellow!30!green] (90:1.2em) to (270:1.2em) to[bend right=60, looseness=1.6] (234:1.2em) to[bend right=60, looseness=1.6] (198:1.2em) (90:1.2em) to[bend left=60] (135:1.2em) (270:1.2em) to[bend left=25] (0:1.2em) to [bend left=25] (90:1.2em);
            \node[circle,fill=white,draw=black,inner sep=0,outer sep=0,minimum size=.25em] at (0:1.2em) {};\node[circle,fill=white,draw=black,inner sep=0,outer sep=0,minimum size=.25em] at (90:1.2em) {};
            \node[circle,fill=white,draw=black,inner sep=0,outer sep=0,minimum size=.25em] at (135:1.2em) {};
            \node[circle,fill=white,draw=black,inner sep=0,outer sep=0,minimum size=.25em] at (198:1.2em) {};
            \node[circle,fill=white,draw=black,inner sep=0,outer sep=0,minimum size=.25em] at (234:1.2em) {};
            \node[circle,fill=white,draw=black,inner sep=0,outer sep=0,minimum size=.25em] at (270:1.2em) {};
            \node[font=\tiny] at (0:.4em) {$\times$};
        \end{scope}
        \end{scope} % end middle
        \begin{scope}[xshift=16em] % right
        \draw[-stealth] (-3.5em,1em) to (-1.5em,2.2em);
        \draw[-stealth] (-3.5em,-1em) to (-1.5em,-2.2em);
        \draw[-stealth] (1.5em,2.2em) to (3.5em,1em);
        \draw[-stealth] (1.5em,-2.2em) to (3.5em,-1em);
        \begin{scope}[xshift=-5em] % L
            \draw[draw=blue, fill=blue] (0,0) circle(.1em);
        \end{scope}
        \begin{scope}[yshift=2em] % T
            \draw[draw=blue, fill=blue] (0,0) circle(.1em);
            \node[circle, minimum size=.1em, outer sep=2pt, inner sep=0] (C) at (0,0) {};
            \node[circle, outer sep=1pt, inner sep=0] (TL) at (120:1.8em) {};
            \node[circle, outer sep=1pt, inner sep=0] (TR) at (60:1.8em) {};
            \draw[-stealth] (C) to (TL);
            \draw[-stealth] (TL) to (TR);
            \draw[-stealth] (TR) to (C);
            %\node at (90:.9em) {$\mu_6$};
            \draw[line cap=round, dash pattern=on 0pt off 1.2pt, line width=.6pt] (112:1.1em) to[bend right=30] (95:1.5em) (68:1.1em) to[bend left=30] (85:1.5em) (110:.7em) to[bend left=30] (68:.7em);
            \node[font=\scriptsize] at (30:1.4em) {$\mu_6$};
        \end{scope}
        \begin{scope}[yshift=-2em] % B
            \draw[draw=blue, fill=blue] (0,0) circle(.1em);
            \node[circle, minimum size=.1em, outer sep=2pt, inner sep=0] (C) at (0,0) {};\node[circle, outer sep=1pt, inner sep=0] (BL) at (240:1.8em) {};
            \node[circle, outer sep=.2pt, inner sep=0] (BR) at (300:1.2em) {};
            \node[circle, outer sep=1pt, inner sep=0] (BRR) at (300:2.4em) {};
            \draw[-stealth] (BRR) to (BR);
            \draw[-stealth] (BR) to (C);
            \draw[-stealth] (C) to (BL);
            \draw[line cap=round, dash pattern=on 0pt off 1.2pt, line width=.6pt] (-110:.7em) to[bend right=30] (-68:.7em) (-68:.9em) to[bend right=60, looseness=1.2] (-63:1.8em);
        \end{scope}
        \begin{scope}[xshift=5em] % R
            \draw[draw=blue, fill=blue] (0,0) circle(.1em);
            \node[circle, minimum size=.1em, outer sep=2pt, inner sep=0] (C) at (0,0) {};
            \node[circle, outer sep=1pt, inner sep=0] (TL) at (120:1.8em) {};
            \node[circle, outer sep=1pt, inner sep=0] (TR) at (60:1.8em) {};
            \draw[-stealth] (C) to (TL);
            \draw[-stealth] (TL) to (TR);
            \draw[-stealth] (TR) to (C);
            \node[font=\scriptsize] at (30:1.4em) {$\mu_6$};
            \draw[line cap=round, dash pattern=on 0pt off 1.2pt, line width=.6pt] (112:1.1em) to[bend right=30] (95:1.5em) (68:1.1em) to[bend left=30] (85:1.5em) (110:.7em) to[bend left=30] (68:.7em);
            \node[circle, outer sep=1pt, inner sep=0] (BL) at (240:1.8em) {};
            \node[circle, outer sep=.2pt, inner sep=0] (BR) at (300:1.2em) {};
            \node[circle, outer sep=1pt, inner sep=0] (BRR) at (300:2.4em) {};
            \draw[-stealth] (BRR) to (BR);
            \draw[-stealth] (BR) to (C);
            \draw[-stealth] (C) to (BL);
            \draw[line cap=round, dash pattern=on 0pt off 1.2pt, line width=.6pt] (-110:.7em) to[bend right=30] (-68:.7em) (-68:.9em) to[bend right=60, looseness=1.2] (-63:1.8em);
        \end{scope}    
        \end{scope} % end right
    \end{tikzpicture}
    %\includegraphics[width=0.5\linewidth]{}
   % \caption{Caption}
    \caption{ Example of a (homotopy) pushout where the gluing is along edges in the (dual) ribbon graph $\zD^*$ on the left, resulting in the gluing of orbifold surfaces in the middle and the gluing of algebras along a vertex (in the homotopy category of dg categories) on the right. The nontrivial $\mu_6$ records the orbifold structure and setting $\mu_6 \equiv 0 $ corresponds to the cosheaf construction of smooth surfaces in Section~\ref{sec: Fukaya}.}
     \label{fig: Cosheaf Orbifold}
\end{figure}

The observations on partial compactifications of surfaces by orbifold points coming from the  perspective of $A_\infty$-deformations of Fukaya categories of surfaces can be developed into an intrinsic theory of partially wrapped Fukaya categories of graded orbifold surfaces $(S, M, \eta)$ with orbifold points of order 2. This has been done independently from various perspectives in \cite{AmiotPlamondon,BarmeierSchrollWang,ChoKim}. Generalizing the notion of dissections to orbifold surfaces, one can define the partially wrapped Fukaya category of a graded orbifold surface in several different ways
\begin{itemize}

\item via explicit $A_\infty$-categories associated to arc systems on orbifold surfaces using tagged arcs in \cite{ChoKim} and untagged arcs in \cite{BarmeierSchrollWang},

\item as global quotient or  invariant categories:  as the orbit 
category $\mathcal W (\widetilde S, \widetilde M, \widetilde \eta) / \mathbb Z_2$  of the Fukaya category of a double cover  
$(\widetilde S, \widetilde M, \widetilde \eta)$ of the orbifold surface in \cite{AmiotPlamondon, BarmeierSchrollWang}
and as the invariant category $\mathcal W (\widetilde S, \widetilde M, \widetilde \eta)^{\mathbb Z_2}$ of a double cover in \cite{ChoKim} (note that $\mathcal W (\widetilde S, \widetilde M, \widetilde \eta) / \mathbb Z_2 \simeq \mathcal W (\widetilde S, \widetilde M, \widetilde \eta)^{\mathbb Z_2}$),

\item as categories of global sections of cosheaves of $A_\infty$-categories of radical square zero dg algebras of types $\mathbb{A}$, $\widetilde{\mathbb{A}}$, and $\mathbb{D}$ glued along a ribbon graph of the orbifold surface in  \cite{BarmeierSchrollWang}.
\end{itemize}
In \cite{BarmeierSchrollWang} it is shown that these three constructions are equivalent. 
Furthermore, the third perspective via cosheaves of $A_\infty$-categories can be viewed as an instance of a conjecture of Kontsevich \cite{Kontsevich}, generalizing the case of smooth surfaces (see Section \ref{sec: Fukaya}).
 In Figure~\ref{fig: Cosheaf Orbifold}, an example of the construction of  $\calw (S,M,\eta)$ as a homotopy pushout of dg categories    along the (red) ribbon graph on the left is given. The interior 3-gon containing an orbifold point gives rise on the right side to higher products $\mu_6$ which can also be viewed as a component of type $\mathbb D$, see \cite{BarmeierSchrollWang,BarmeierWangICRA}.

\subsection{New graded associative algebras from Fukaya categories of orbifold surfaces}

The Fukaya category $\mathcal W (S, M, \eta)$ of a graded orbifold surface always admits a formal generator whose corresponding associative algebra is graded skew-gentle. This gives a symplectic interpretation of the different orbifold surface models for the derived categories of  skew-gentle algebras \cite{AmiotDerivedSkewGentle,AmiotBruestle, LabardiniSchrollValdivieso, QiuZhangZhou}. Unlike for gentle algebras, the class of skew-gentle algebras is not closed under derived equivalence. However, the algebras derived equivalent to  type $\mathbb D$ quivers (which are derived skew-gentle) were classified  in \cite{AmiotPlamondon} using the geometric model of orbifold disks. More generally, dissections of orbifold surfaces give rise to a new class of associative graded algebras called \textit{semi-gentle algebras}. These algebras are classified in terms of orbifold dissections in \cite{BarmeierSchrollWang,KimPhD}. Semi-gentle algebras are derived equivalent to graded skew-gentle algebras and are conjectured to be closed under derived equivalence \cite[Conjecture 8.1]{BarmeierSchrollWang}.

\section{Outlook}

Geometric models for graded gentle and skew-gentle algebras give a comprehensive understanding of their derived categories and emphasise the importance of considering graded algebras in representation theory.   The derived tame (graded) hereditary algebras of types $\mathbb A$ and $\mathbb D$ give rise to (graded) gentle and skew-gentle algebras via gluings of polygons or orbifold polygons (realised as homotopy colimits of dg categories) and this construction covers most known derived tame algebras. 
An  exception are the algebras of type $\mathbb E$. They are hereditary and representation-finite, hence derived discrete and to the best of our knowledge, no symplectic surface model currently exists.

Moreover, it is not known what conditions are needed, for example, to ensure that a homotopy colimit of derived tame algebras remains derived tame. (Co)sheaf-theoretical techniques have a solid foundation in the homotopy theory of dg  and $A_\infty$-categories \cite{KellerDG,Tabuada,Toen} as well as $\infty$-categories \cite{Lurie}, but their importance for representation theory is only starting to be understood. The mounting evidence such as \cite{BarmeierSchrollWangDef,BarmeierSchrollWang,ChristHaidenQiu,DyckerhoffKapranov18,HaidenKatzarkovKontsevich} as well as recent advances in cluster theory \cite{Christ,KellerLiu} showcase the importance of these techniques. 

Whereas ``most'' derived tame algebras have a symplectic surface model, corresponding via homological mirror symmetry to nodal curves in the sense of \cite{BurbanDrozdnodal,LekiliPolishchuk}, derived tameness cannot necessarily be expected to be found in higher dimensional geometry where it is a reasonable expectation that most algebras are derived wild. However, this need not preclude the effectiveness of symplectic methods in representation theory. Associative algebras arise in particular from formal generators in triangulated categories and their representation theory and $A_\infty$-deformations thus can be expected to appear in a wide variety of contexts. Examples  of geometric models for the derived categories of associative algebras appear for example in \cite{DyckerhoffJassoLekili}, where it is shown that (derived) partially wrapped Fukaya categories of symmetric powers of disks with stops are equivalent to perfect derived categories of higher Auslander algebras of type $\mathbb A$. This relates the symplectic geometry of higher-dimensional symplectic manifolds to higher-dimensional Auslander--Reiten theory in representation theory \cite{Iyama}. Another higher-dimensional example is given by extended Khovanov arc algebras, whose presentation by quivers and relations is given in \cite{BarmeierWangKho}, see also \cite{Bowman1}. In \cite{MakSmith} it is shown that the perfect derived categories of these graded algebras are equivalent to the derived Fukaya--Seidel categories of Hilbert schemes of Milnor fibers. 

It can be hoped that the interactions between representation theory and symplectic geometry advance our understanding both of higher dimensional geometry as well as that of the representation theory of wild algebras and their derived categories.

\bibliographystyle{acm}
\bibliography{ICMsurvey}
\end{document}